\documentclass{amsart}

\usepackage{amsthm,amsfonts,amsmath,amssymb,latexsym,epsfig}

\newtheorem{theorem}{Theorem}
\newtheorem{lemma}[theorem]{Lemma}
\newtheorem{proposition}[theorem]{Proposition}

\newenvironment{remark}{\medskip \refstepcounter{theorem}
\noindent  {\bf Remark \thetheorem}.\rm}{\,}
\newenvironment{definition}{\medskip \refstepcounter{theorem}
\noindent  {\bf Definition \thetheorem}.\rm}{\,}

\renewcommand{\thetheorem}{\thesection.\arabic{theorem}}

\theoremstyle{definition}

\theoremstyle{remark}

\def \ZN{{\mathbb Z}[1/N,\zeta_N]}
\def \o{\omega}

\def\wh{\widehat{\phantom{p}}}

\def \tf{\tilde{f}}

\def \bZ{{\mathbb Z}}
\def \bG{{\mathbb G}}

\def \D{\{\d,\dd\}}
\def \dd{\delta_q}

\def \<{\langle}
\def \>{\rangle}

\def \cU{\mathcal U}

\def \d{\delta_p}

\def \bZ{{\mathbb Z}}

\def \bC{{\mathbb C}}

\def \ra{\rightarrow}

\def \bq{{\bf q}}

\def \bQ{{\mathbb Q}}

\def \bQ{{\mathbb Q}}

\def \bM{{\mathbb M}}

\def \bZ{{\mathbb Z}}
\def \bG{{\mathbb G}}

\begin{document}

\title[Arithmetic PDES on modular curves]{Arithmetic partial differential
equations, II:  modular curves}
\author[A. Buium]{Alexandru Buium}
\thanks{During the preparation of this work, the first author
was partially supported by NSF grant DMS 0552314.}
\author[S.R. Simanca]{Santiago R. Simanca}
\address{Department of Mathematics and Statistics \\
University of New Mexico \\ Albuquerque, NM 87131}
\email{buium@math.unm.edu, santiago@math.unm.edu}

\begin{abstract}
We continue the study of arithmetic partial differential equations
initiated in \cite{pde} by classifying ``arithmetic convection equations'' on
modular curves, and by describing their space of solutions. Certain of these
solutions involve the Fourier expansions of the Eisenstein modular
forms of weight $4$ and $6$, while others involve the
Serre-Tate expansions \cite{Mori, shimura} of the same modular
forms; in this sense, our arithmetic convection equations can be seen as
``unifying" the two types
of expansions. The theory can be generalized to one of
 ``arithmetic heat equations'' on modular curves, but we prove that
modular curves do not carry ``arithmetic wave equations.'' Finally, we prove
an instability result for families of arithmetic heat equations converging to
an arithmetic convection equation.
\end{abstract}

\maketitle

\section{Introduction}
\subsection{Main concepts and results}
This article is a direct continuation of \cite{pde}, and for the purpose
of the current introduction, we shall assume some familiarity with the
introduction to our earlier work. In there 
we developed an arithmetic analogue of the theory
of partial differential equations in space-time, where the r\^ole of the
spatial derivative is played by a Fermat quotient
operator $\d$ with respect to a fixed prime $p$, while the r\^ole of
derivative with respect to (the ``exponential'' $q$ of) time is
played by the usual derivation operator $\dd$. (We view the
operator $\d$ as a derivative with respect to $p$; the ordinary
differential theory based on this interpretation was initiated in
\cite{char}.) The rings on which $\d$ and $\dd$ act, which satisfy
certain compatibility conditions, are called $\D$-rings. If
$X$ is a scheme of finite type over a $p$-adically complete
$\D$-ring $A$, we define a partial differential operator as a map
$f:X(A) \ra A$ given locally, in the Zariski topology, by a restricted
$p$-adic power series in the affine coordinates $x$ and finitely
many of its iterated ``derivatives'' $\d^i\dd^j x$.

Given a partial differential operator $f:X(A) \ra A$, we may
consider the {\it equation} $f=0$, and its {\it space of solutions}
$\cU:=f^{-1}(0) \subset X(A)$.
 It is natural to restrict our
attention to special cases where an appropriate notion of {\it
linearity} for $f$ can be defined, and then to investigate such
linear operators. In our earlier work \cite{pde} we concentrated primarily
on the case when $X$ is a commutative group scheme of dimension $1$
(especially ${\mathbb G}_a$, ${\mathbb G}_m$, or $X=E$, an elliptic
curve) and $f$ is a homomorphism (to the additive group of $A$); the
homomorphism condition was meant to be an analogue of linearity for
classical partial differential operators in analysis.
Three classes of these equations ended up playing a
central r\^ole in our study. They can
be viewed as arithmetic analogues of convection, heat, and wave equations,
respectively.

In the current article we consider a modular analogue of the situation above,
where we now take $X$ to be  the total space $\bM_{\Gamma_1(N)}$ of
a natural ${\mathbb G}_m$-torsor over a modular curve, or a
``stacky'' version of such an object. (With precision,
$\bM_{\Gamma_1(N)}$ classifies elliptic curves, possibly with some level
structure, together with an invertible $1$-form. In our work here,
we present the theory without reference to
 modular curves, in a purely stacky manner, in the style of
\cite{Katz} and  \cite{difmod}; the translation in terms of modular
curves is straightforward.) We shall assume that the function $f$ is
homogeneous of a certain weight with respect to the action of the
multiplicative group on $X$; this will make of $f$ a {\it partial
differential modular form} (cf. \cite{difmod} for the ordinary
differential case). For such an $f$, the space of solutions $\cU$ has
a ${\mathbb G}_m$-action. In addition, we will assume that $f$ has a
certain covariance with respect to the action of Hecke
correspondences; these $f$s will be called  {\it isogeny covariants}
(again, cf. \cite{difmod} for the ordinary differential case), and
for them, the space of solutions $\cU$ is also saturated with
respect to ``isogeny.'' Isogeny covariant forms of {\it covariance
degree} $1$ will be viewed as arithmetic analogues, in our context
here, of  linear partial differential equations in analysis. For
forms of even integral weight $m$, the covariance degree equals
$-m/2$, and so linearity corresponds to weight $-2$.
In what follows, we will consider a generalization of this situation,
where the $f$s are allowed to have singularities along
the supersingular locus ``$E_{p-1}=0$" of $X$.

In parallel with our results in \cite{pde}, we address two main
problems here: the first is to find all possible linear partial
differential equations, and the second is to describe the space of
solution $\cU$ of any one such. We deal with these two problems in 
sections \S 2 and \S 3, respectively.
The ordinary differential
theory \cite{difmod} in the arithmetic $p$-direction (respectively,
in the geometric $q$-direction) provides
 order $1$ isogeny covariant forms $f^1_p$ and $f^{\partial}_p$
 involving $\d$ only (respectively, an order $1$
isogeny covariant form $f^1_q$ involving $\dd$ only). Both $f^1_q$
and the product $f^{1}_{p,-2}:=
f^{\partial}_p f^1_p$ have weight $-2$.  In \S 2 of
we prove that $f^1_q$ and $f^1_{p,-2}$  form
a basis of the space of isogeny covariant forms with
singularities along $E_{p-1}=0$, of weight $-2$, and order $1$.
These forms should be viewed as ``arithmetic convection equations.''
In \S 3 we analyze the spaces of solutions $\cU$ for  such forms. We
are particularly interested in two subsets of $\cU$, denoted by
$\cU_{bad}$ and $\cU_{good}$, respectively, that consist of solutions with
bad and good reductions at infinity, and which satisfy certain
non-degeneracy conditions.  If  $f$ is a linear combination of the form
$f=f^1_q+\lambda f^1_{p,-2}$, then, for $\lambda$
invertible, we have the following ``quantization'' phenomenon
(reminsicent of our theory in \cite{pde}): if the parameter
$\lambda$ is not a positive integer times a certain fixed quantity,
then the set $\cU_{bad}$ is a union of ${\mathbb G}_m$-orbits of one
basic solution (that can be expressed in terms of the Fourier expansions of
Eisenstein modular forms of
weights $4$ and $6$) and of their shifts by roots of unity. For the
exceptional values of $\lambda$, the space $\cU_{bad}$ has a
different, more complicated structure, which we shall completely
describe in terms of closed form formulae.
With the set
$\cU_{good}$ we encounter a similar phenomenon: for $\lambda$ not
a positive integer times a certain fixed quantity,
$\cU_{good}$ consists of ${\mathbb G}_m$-orbits of ``stationary''
solutions. For the exceptional values of $\lambda$, the set
$\cU_{good}$ has again a richer but completely understood
structure. The r\^ole of the Fourier series of Eisenstein forms of
weights $4$ and $6$
is played now, in the good reduction case, by the ``Serre-Tate
expansions'' \cite{Mori, shimura} of these modular forms. In this way,
both the Fourier expansions and the Serre-Tate expansions appear in
the solutions of the same arithmetic partial differential equation;
in a sense, these two types of expansions are ``unified'' by the said
equation. Furthermore, we will discover ``canonical'' $1$-to-$1$
correspondences between certain sets of bad reduction
solutions and certain sets of good reduction
solutions of appropriate pairs of arithmetic convection equations.
A similar result holds for ``arithmetic heat equations'' that
are ``close to'' arithmetic convection equations.

Although arithmetic convection and heat equations exist,
we will prove also that, in a suitable sense, modular curves do not carry
``arithmetic wave equations.'' 

We shall end with the discovery of an ``instability phenomenon.''
As  ``arithmetic heat equations'' converge to a given
``arithmetic convection equation,''
the solutions of the arithmetic heat equations with given
boundary conditions {\it do not converge} to the corresponding
solutions of the arithmetic convection equation. This
situation does not seem to have a parallel when viewing
solutions of heat versus convection equations
in real analysis.

\subsection{Review of terminology and notation}
 We end this introduction by recalling some of the basic concepts and
notation from \cite{pde} that we will need later on, as well as introducing
the ring of weights. 

A derivation from a
ring $A$ to an $A$-algebra $B$ is a map $\delta:A \ra B$ such that
$\delta(x+y)=\delta x + \delta y$ and $\delta(xy)=x \delta y + y
\delta x$. If $p$ is a prime, a $p$-derivation is a map $\delta=\delta_p : A
\ra B$ such that
$$\begin{array}{rcl}
\delta(x+y) & = &  \delta x + \delta y +C_p(x,y)\, , \\
\delta(xy)  & = & x^p  \delta y +y^p  \delta x  +p  \delta x  \delta
y\, ,
\end{array}$$
where $C_p(X,Y)$ stands for the polynomial with ${\mathbb Z}$-coefficients
$$C_p(X,Y):=\frac{X^p+Y^p-(X+Y)^p}{p}\, .$$
If $\d$ is a $p$-derivation
then $\phi_p:A  \ra   B $ defined by
$$
\phi_p x  :=  x^p+p\d x\, ,
$$
is always a ring
homomorphism. A $\D$-ring is a ring $A$ equipped with a derivation
$\dd:A \ra A$ and a $p$-derivation $\d:A \ra A$, such that
\begin{equation}
\label{compp}
\dd \d x=p\d \dd x + (\dd x)^p-x^{p-1}\dd x\, .
\end{equation}
A morphism of $\D$-rings is a ring homomorphism that commutes with
$\d$ and $\dd$.
 A $\D$-{\it prolongation sequence}
is a sequence of rings  $S^{*}=\{S^n\}_{n\geq 0}$ equipped with ring
homomorphisms $\varphi: S^n \ra S^{n+1}$ (used to view each
$S^{n+1}$ as an $S^n$-algebra), $p$-derivations $\d:S^n \ra
S^{n+1}$, and derivations $\dd :S^n \ra S^{n+1}$ satisfying
condition (\ref{compp}), and such that $\d \circ \varphi=\varphi \circ
\d$, and $\dd \circ \varphi=\varphi \circ \dd$, respectively. A {\it morphism}
$u^*:S^* \ra \tilde{S}^*$ of $\D$-prolongation sequences is a
sequence $u^n:S^n \ra \tilde{S}^n$ of ring homomorphisms such that
$\d \circ u^n=u^{n+1} \circ \d$, $\dd \circ u^n=u^{n+1} \circ \dd$,
and $\varphi \circ u^n=u^{n+1} \circ \varphi$. If $A$ is a $\D$-ring,
then we can attach to it a prolongation sequence (which we still
denote by $A$) by taking $A^n=A$, $\varphi=id$. A prolongation
sequence $S^*$ over a $\D$-ring $A$ is a morphism of
$\D$-prolongation sequences $A \ra S^*$.

Throughout our work we fix a prime $p\geq 5$. A basic example of
$\D$-ring is  $A=R:=R_p=({\mathbb Z}_p^{ur})\wh$, with $\dd=0$,
and
$$\d x=\frac{\phi_p(x)-x^p}{p}\, ,$$
where $\phi_p$ is the unique lift of the Frobenius map to
$R/pR$. Here, and in the sequel, the superscript $\wh$ denotes
$p$-adic completion. Another basic example of $\D$-ring that plays a
r\^ole here is $A=R((q))\wh$, where $q$ is an indeterminate over $R$,
$R((q)):=R[[q]][q^{-1}]$, $\dd F=q \frac{\partial F}{\partial q}$,
and
$$\d F=\frac{F^{(\phi_p)}(q^p)-(F(q))^q}{p}\, ,
$$
where the superscript $(\phi_p)$
means twisting coefficients by $\phi_p:R \ra R$. A basic example of
prolongation sequence is as follows.  For a $\D$-ring $A$ and $z$ a
variable, we define
$$A[z^{(\leq r)}]\wh:=A[z^{(i,j)}|_{i+j \leq r}]\wh$$ where $i,j
\geq 0$, $z^{(0,0)}=z$. Then there is a unique structure of a
$\D$-prolongation sequence on this sequence, extending that of $A$,
such that $\d z^{(i,j)}=z^{(i+1,j)}$ and $\dd
z^{(i,j)}=z^{(i,j+1)}$. A similar definition can be given for $z$ a
tuple of variables.  More generally, if $g \in A[z]\wh\backslash
(p)$, we may consider the sequence  $A[z^{(\leq r)},g^{-1}]\wh$; this
has a unique structure of $\D$-prolongation sequence extending the
one just defined. 

We will also need the ring of {\it weights}
$W:={\mathbb Z}[\phi]$, the polynomial ring in the ``variable''
$\phi$; we let $W_+$ denote the semigroup of all $w=\sum_{i=0}^r a_i
\phi^i \in W$ with $a_i \geq 0$. We define ${\rm ord}(w)=r$ if $a_r
\neq 0$, $deg(w)=\sum a_i$. If $\lambda$ is an invertible element in
a $\D$-ring $A$ and $w \in W$ (or if $\lambda$ is not necessarily
invertible but $w \in W_+$), we set $\lambda^w=\prod
\phi_p^i(\lambda)^{a_i}$. The same notation can be introduced for
$\D$-prolongation sequences.

If instead of $\D$ we use only $\d$ (respectively $\dd$) we are led
to obvious notions of $\d$-ring, $\d$-prolongation sequence, etc.
(respectively $\dd$-ring, $\dd$-prolongation sequence, etc.).

\medskip

\noindent {\bf Acknowledgements.}  We are grateful to Dinesh Thakur
for discussions on the Artin-Hasse and quantum exponentials.

\section{Equations}
\setcounter{theorem}{0}
In this section we introduce and classify some of our arithmetic
analogues of linear partial differential equations on modular
curves.

\subsection{$\D$-modular forms}

The following definitions are analogues in the partial differential
case of some definitions given in \cite{difmod}; they follow the
viewpoint put forward in \cite{Katz} that treated the
``non-differential'' setting.

Let $N$ be a positive integer that is not divisible by $p$, and that
we shall fix throughout the paper. Given any ring $S$ where $6$ is an
invertible
element, we denote by $\bM(\Gamma_1(N),S)$ the set of all triples
$(E/S,\alpha,\omega)$ where $E/S$ is an elliptic curve, $\omega$ is an
invertible $1$-form on $E$, and $\alpha \, : \, ({\mathbb Z}/N{\mathbb Z})_S
\ra E$ is a closed immersion of group schemes (referred to as a
$\Gamma_1(N)$-{\it level structure}). When $N=1$, we shall usually drop
$\alpha$ and $\Gamma_1(1)$ from the notation. In particular, we write
$\bM(S)$ for $\bM(\Gamma_1(1),S)$.

Notice that we have an identification
$$\bM(S)=\{(a,b)\in S \times S\ |\ 4a^3+27b^2 \in S^{\times}\}$$
given by attaching to each pair $(a,b) \in \bM(S)$, the pair
$(E,\omega)$ where $E$ is the projective closure of the affine curve
$y^2=x^3+ax+b$, and $\omega=dx/y$.

There is an action of ${\mathbb G}_m(S)=S^{\times}$ on
$\bM(\Gamma_1(N),S)$ via multiplication on the $\omega$ component.
The induced action on $\bM(S)$ is given by the formula
$\lambda(a,b)=(\lambda^{-4} a, \lambda^{-6} b)$.

\begin{remark}
Notice that when $N>4$, $\bM(\Gamma_1(N),S)$ is the set of
$S$-points of a natural ${\mathbb G}_m$-torsor $\bM_{\Gamma_1(N)}$
over the modular curve $Y_1(N)$. The same is true when $N \leq 4$,
but only in a ``stacky sense.'' Cf. \cite{Katz}, \cite{DI}. For
$N=1$, $\bM(S)$ is, of course, the set of $S$-points of the
{\it modular scheme}
$$\bM=\bM_{\Gamma}:={\rm Spec}\
\bZ[1/6][a_4,a_6,(4a_4^3+27a_6^2)^{-1}],$$ where $a_4,a_6$ are
indeterminates and $\Gamma$ stands for ${\rm SL}_2(\bZ)=\Gamma_1(1)$.
\qed
\end{remark}
\medskip

Let us fix a weight $w$ in the ring of weights $W$, with ${\rm ord}(w) \leq r$,
and a $\D$-ring $A$ that is Noetherian, $p$-adically complete, and an
integral domain of characteristic zero.

\begin{definition}
A $\D$-{\it modular form} over  $A$, of weight $w \in W$
and order $r$ on $\Gamma_1(N)$, is a rule $f$
  that associates to any $\D$-prolongation sequence $S^*$
  of Noetherian, $p$-adically complete rings over $A$,
  and any triple $(E/S^0, \alpha, \o) \in \bM(\Gamma_1(N),S^0)$,
  an element $f(E/S^0, \alpha, \o, S^*) \in S^r$ subject to the following
  conditions:
\begin{enumerate}
 \item $f(E/S^0, \alpha, \o, S^*)$ depends on the isomorphism class
  of  $(E/S^0, \alpha, \o)$ only.
\item  The formation of  $f(E/S^0, \alpha, \o, S^*)$  commutes with base change
$u^*:S^* \ra \tilde{S}^*$, that is to say,
$$f(E \otimes_{S^0}
 \tilde{S}^0/\tilde{S}^0, \alpha \otimes \tilde{S}^0,
 u^{0*} \o, \tilde{S}^*)=
 u^r(f(E/S^0, \alpha, \o, S^*))\, .$$
\item $f(E/S^0, \alpha, \lambda \o, S^*)=\lambda^{-w} \cdot
 f(E/S^0, \alpha, \o, S^*)$
for all $\lambda \in (S^0)^{\times}$. \qed
\end{enumerate}
\end{definition}

We denote by $M^r_{pq}(\Gamma_1(N),A,w)$ the $A$-module of all
$\D$-modular forms over $A$ of weight $w \in W$ and order $r$ on
$\Gamma_1(N)$.

We have natural maps
$$\varphi: {M}^{r}_{pq}(\Gamma_1(N),A,w) \ra
{M}^{r+1}_{pq}(\Gamma_1(N),A,w)$$ that send any $f$ into $\varphi
\circ f$, and we have natural maps
$$\phi: {M}^{r}_{pq}(\Gamma_1(N),A,w) \ra {M}^{r+1}_{pq}(
\Gamma_1(N),A,\phi w)$$ that send any $f$ into $f^{\phi}:=\phi \circ
f$. Finally, notice that for any integer $N^{'}$ such $N'|N$, there are
natural maps
$${M}^{r}_{pq}(\Gamma_1(N'),A,w)\ra {M}^{r}_{pq}(
\Gamma_1(N),A,w)\, .$$
In particular, for $N'=1$, we have a map
$${M}^r_{pq}(A,w) \ra {M}^{r}_{pq}(\Gamma_1(N),A,w).
$$

Following \cite{difmod} verbatim, we obtain the following description of
${M}^r_{pq}(A,w)$. Let $a_4$, $a_6$, and $\Lambda$ be variables, and
set $\Delta:=-2^6a_4^3-2^43^3a_6^2$. Then ${M}^r_{pq}(A,w)$
identifies with the set of all power series $f$ in
$$M^r_{pq}:=A[a_4^{(\leq r)}, a_6^{(\leq r)},\Delta^{-1}]\wh$$
such that
$$\begin{array}{r}
f(\ldots ,\d^i \dd^j (\Lambda^4 a_4), \ldots ,\d^i
 \dd^j (\Lambda^6 a_6),\ldots ,
\Lambda^{-12}\Delta^{-1}) = \hspace{1.5in} \mbox{} \\
\Lambda^w f(\ldots , \d^i \dd^j  a_4, \ldots ,\d^i \dd^j  a_6,
\ldots , \Delta^{-1})\, .
\end{array}
$$
Given $f \in {M}^r_{pq}(A,w)$, we shall use the same symbol and
still denote by $f\in M^r_{pq}$ the corresponding series.
The map
\begin{equation} \label{othf}
\bM(A) \ra A
\end{equation}
defined by
$$
(a,b)\mapsto f( \ldots ,\d^i \dd^j a, \ldots ,\d^i \dd^j b, \ldots ,
\Delta^{-1})\, ,
$$ will be denoted by $f$; it has  the property
that
$$
f(\lambda^4 a,\lambda^6 b)=\lambda^w f(a,b)
$$
for all $\lambda \in A^{\times}$.

A form  $f \in {M}^r_{pq}(\Gamma_1(N),A,w)$ will be called {\it
essentially of level one} if the rule $f$ does not depend on the
variable $\alpha$. We denote by  $M^r_{pq}(\Gamma_1(N),A,w)_1$ the
subspace of ${M}^r_{pq}(\Gamma_1(N),A,w)$ consisting of forms that
are essentially of level one.

\begin{lemma}
\label{allfa} The natural maps ${M}^r_{pq}(A,w) \ra
{M}^r_{pq}(\Gamma_1(N),A,w)_1$ are isomorphisms.
\end{lemma}

{\it Proof}. We prove surjectivity. Injectivity follows
similarly.

Let $$\tf \in {M}^r_{pq}(\Gamma_1(N),A,w)_1\, ,$$ and consider a
triple $(E,\omega,S^*)$. There is a Galois \'{e}tale $S^0$-algebra
$\tilde{S}^0$ such that $E \otimes_{S^0} \tilde{S}^0$ has a
$\Gamma_1(N)$-level structure $\alpha$. As in (3.15) of \cite{difmod},
$S^* \otimes_{S^0}\tilde{S}^0$ has a unique structure of
$\D$-prolongation sequence extending that of $S^*$. Set
$$
f(E,\omega,S^*):=\tf(E \otimes_{S^0}\tilde{S}^0, \omega \otimes 1,
\alpha, S^* \otimes_{S^0}\tilde{S}^0)\, .
$$
This quantity does not
depend on the choice of $\alpha$, and by functoriality is seen to
belong to $S^r$ (rather than $S^r
\otimes_{S^0}\tilde{S}^0$). Clearly $f$ is mapped into $\tf$.
 \qed
\medskip

Now for any $f \in {M}^r_{pq}(A,w)$, we may consider the map
$f:\bM(A) \ra A$ in (\ref{othf}). We say that $f=0$ is an
equation of weight $w$, and $f^{-1}(0) \subset \bM(A)$ is its {\it
space of solutions}. Notice that $f^{-1}(0)$ is invariant under the
${\mathbb G}_m$-action:  if $(a,b)\in f^{-1}(0)$ and $v \in
A^{\times}$, then $(v^4 a, v^6 b)\in f^{-1}(0)$. In what follows, we
will  impose another property on $f$, {\it isogeny covariance},
which gives the space of solutions an extra symmetry property with
respect to ``Hecke correspondences.''

\subsection{Isogeny covariance} Let $f \in
{M}^r_{pq}(\Gamma_1(N),A,w)$ be a $\D$-modular form of weight
$w=\sum n_i \phi^i$ on $\Gamma_1(N)$. Assume ${\rm deg}(w):=\sum
n_i$ is even. Generalizing the level one definition in
\cite{difmod}, we say that $f$ is {\it isogeny covariant} if for any
$\D$-prolongation sequence $S^*$, any triples
$(E_1,\alpha_1,\omega_1),(E_2,\alpha_2,\omega_2) \in
\bM(\Gamma_1(N),S^0)$, and
 any isogeny $u:E_1 \ra E_2$ of degree prime to $p$,
 with $\omega_1=u^* \omega_2$ and $u \circ \alpha_1=\alpha_2$, we have that
$$f(E_1,\alpha_1,\omega_1,S^*)={\rm deg}(u)^{-{\rm deg}(w)/2} \cdot
f(E_2,\alpha_2, \omega_2,S^*)\, .$$
The number $-{\rm deg}(w)/2$ will be
called the {\it covariance degree} of $f$. Forms of covariance
degree $1$ will be called {\it linear}; they should be viewed as
the analogues, in our current context, of the linear differential operators in
analysis.

We denote by $I^r_{pq}(\Gamma_1(N),A,w)$ the space of isogeny
covariant forms belonging to the space ${M}^r_{pq}
(\Gamma_1(N),A,w)$. We denote by  $I^r_{pq}(\Gamma_1(N),A,w)_1$ the
subspace of $I^r_{pq}(\Gamma_1(N),A,w)$  consisting of forms which
are essentially of level one. By Lemma \ref{allfa}, it follows that
the natural maps $I^r_{pq}(A,w) \ra I^r_{pq}(\Gamma_1(N),A,w)_1$ are
isomorphisms.

Notice that if $f$ is isogeny covariant, then its space of solutions
$f^{-1}(0)$ has the following extra symmetry (which is morally a
symmetry with respect to Hecke correspondences). For if we assume
$y^2=x^3+a_1x+b_1$ and $y^2=x^3+a_2x+b_2$ are two elliptic curves
with coefficients in $A$, and that there exists an isogeny over $A$
of degree prime to $p$ between them that pulls back $dx/y$ into
$dx/y$, then $(a_1,b_1) \in f^{-1}(0)$ if, and only if, $(a_2,b_2)
\in f^{-1}(0)$.

\subsection{Variants}
If in the definitions above we use
$\d$-rings and $\d$-prolongation
sequences instead of $\D$-rings and
$\D$-prolongation sequences,
we obtain the concept of {\it $\d$-modular form}
(over a $\d$-ring $A$); in that case, we denote the corresponding spaces
by ${M}^r_p(\Gamma_1(N),A,w)$ and $I^r_p(\Gamma_1(N),A,w)$,
respectively. For $N=1$, these $\d$-modular forms are the
$\delta$-modular forms in the arithmetic setting of \cite{difmod}.
Similarly, if in these definitions
we use $\dd$-rings and
$\dd$-prolongation sequences {\it of not necessarily $p$-adically
complete rings} instead, we  obtain the concept of {\it $\dd$-modular form}
(over a $\dd$-ring $A$); we then denote the corresponding spaces
by ${M}^r_q(\Gamma_1(N),A,w)$ and $I^r_q(\Gamma_1(N),A,w)$,
respectively. For $N=1$, these $\dd$-modular forms are the
$\delta$-modular forms in the geometric setting of \cite{difmod}.

As before, we can speak of  forms in ${M}^r_p(\Gamma_1(N),A,w)$ and
$I^r_p(\Gamma_1(N),A,w)$ that are essentially of level one. We
denote them by
$${M}^r_p(\Gamma_1(N),A,w)_1\, , \quad I^r_p(\Gamma_1(N), A, w)_1\, ,$$
and once again, these spaces are isomorphic
to ${M}^r_p(A,w)$ and $I^r_p(A,w)$, respectively.
We introduce a similar notation in the case of $\dd$-modular forms.

There are natural injective maps
$$\begin{array}{c}
{M}^r_p(\Gamma_1(N),A,w) \ra {M}^r_{pq}(\Gamma_1(N),A,
w)\, , \\
{M}^r_q(\Gamma_1(N),A,w) \ra {M}^r_{pq}(\Gamma_1(N),A,w)\, .
\end{array}$$
A similarl statement holds for the $I$ spaces in place of the $M$ spaces
above.

 If in all of the definitions above we insist that elliptic
curves have ordinary reduction mod $p$, then we get a new set of
concepts that will be indicated using {\it ord} in the notation. So,
for instance, instead of the sets $\bM(A)$, ${M}^r_{pq}(A,w)$,
$I^r_{pq}(A,w)$, etc., we will have sets $\bM_{ord}(A)$, ${M}^r_{pq,
ord}(A,w)$, $I^r_{pq, ord}(A,w)$, etc. Lemma \ref{allfa} continues
to hold for these ``{\it ord}'' sets. We also have natural maps
$\bM_{ord}(A) \ra \bM(A)$, ${M}^r_{pq}(A,w) \ra
{M}^r_{pq,ord}(A,w)$, $I^r_{pq}(A,w) \ra I^r_{pq,ord}(A,w)$, etc.

Notice that
$$
\bM_{ord}(A)=\{(a,b) \in \bM(A)\ ;\ E_{p-1}(a,b) \in
A^{\times}\}\, ,
$$
where $E_{p-1} \in \bZ_p[a_4,a_6]$ is the polynomial
whose Fourier series is the normalized Eisenstein series
$E_{p-1}(q)$ of weight $p-1$. Here {\it normalized} means that
the constant coefficient is $1$.

\subsection{The main examples}
\label{ex} Let $A$ be a $p$-adically complete
$\D$-ring that is an integral Noetherian
domain of characteristic zero, and let $N \geq 1$. The construction
in \cite{pde}, Remark 8.3, provides a recipe that attaches to any
$\D$-prolongation sequence $S^*$ of Noetherian $p$-adically complete
rings over $A$, and to any pair $(E,\omega)$ consisting of an
elliptic curve $E/S^0$ and an invertible $1$-form $\omega$ on $E$,
elements
$$f^1_p(E,\omega,S^*)\, , f^1_q(E,\omega,S^*) \in S^1\, .$$
Moreover, in the case where $E/S^0$ has ordinary reduction mod $p$, we
define (cf. \cite{barcau}, Construction 3.2 and Theorem 5.1) an
element
$$
f^{\partial}_p(E,\omega,S^*) \in (S^1)^{\times}\, .
$$
We shall sometimes drop $S^*$ from the notation.

The rules $f^1_p$ and $f^1_q$ define isogeny covariant modular forms,
$$
f^1_p \in I^1_p(A,-1-\phi),\quad f^1_q \in I^1_q(A,-2)\, ,
$$
cf. \cite{difmod}, Construction 4.1 and
Construction 4.9. In particular $f^1_p$ and $f^1_q$ give rise to
elements of $I^1_{pq}(A,-1-\phi)$ and $I^1_{pq}(A,-2)$, respectively,
and hence, to elements in $I^1_{pq,ord}(A,-1-\phi)$ and
$I^1_{pq,ord}(A,-2)$, respectively,
 which we continue to denote by $f^1_p$ and $f^1_q$.
Similarly $f^{\partial}_p$ defines an element
$$
f^{\partial}_p \in I^1_{p,ord}(A,\phi-1)\, ,
$$
hence an element of $I^1_{pq,ord}(A,\phi-1)$, which we still denote by
$f^{\partial}_p$.

\begin{remark}
\label{cano}
Let $w_1, w_2 \in W$ be two weights of order $\leq
r$ such that $deg(w_1)=deg(w_2) \in 2 \bZ$. Then $w_2-w_1 = (\phi-1)
w$ for some $w \in W$, and the map
$$
\begin{array}{rcl}
I_{pq,ord}^r(A,w_1) & \ra & I_{pq,ord}^r(A,w_2) \\
g & \mapsto &
(f^{\partial}_p)^w g
\end{array}
$$
is clearly an isomorphism. In particular, we may consider the forms
$$
f^{1\phi^i}_{p,-2}:=(f^{\partial}_p)^{\frac{\phi^{i+1}+\phi^i-2}{\phi-1}}
(f^1_p)^{\phi^i} \in I^s_{p,ord}(A,-2)\, .
$$
Similarly, we may consider the forms
$$
f^{1\phi^i}_{p,-1-\phi^s}=(f^{\partial}_p)^{\frac{\phi^{i+1}+\phi^i-1-\phi^s}
{\phi-1}} (f^1_p)^{\phi^i}\in I^s_{p,ord}(A,-1-\phi^s)\, .
$$
\qed
\end{remark}

\begin{theorem}
\label{gama}
\begin{enumerate}
\item For any $s \geq 1$, the space $I^s_{p,ord}(A,-2) \otimes L$ has an
$L$-basis consisting of the forms
$$
f^1_{p,-2}, f^{1\phi}_{p,-2}, \ldots , f^{1\phi^{s-1}}_{p,-2}\, .
$$
\item The space $I^1_q(A,-2)\otimes L$ has a basis consisting of
$f^1_q$.
\item For any $s \geq 1$, let
$$f^s_p:=\sum_{i=0}^{s-1} p^{s-1-i}f^{1\phi^i}_{p,-1-\phi^s}\in I^s_{p,ord}
(A,-1-\phi^s)\, .$$
 Then $f^s_p$ belongs
to $I^s_p(A,-1-\phi^s)$, and it constitutes a basis for it.
\item For any $s \geq 1$, the space $I^s_{p,ord}(A,-2) \otimes L$ has an
$L$-basis
consisting of the forms
$$
f^1_{p,-2}, f^2_{p,-2}, \ldots , f^{s}_{p,-2}\, .
$$
\end{enumerate}
\end{theorem}

{\it Proof}. (1) is implicit in \cite{barcau}; it also follows
directly by Proposition 8.75 in \cite{book}. (2) is proved in
\cite{difmod}, Corollary 7.24. (3) is implicit in \cite{barcau}, and
also follows by Proposition 8.61 and Theorem 8.83 in \cite{book}.
(4) is a clear consequence of (1) and (3).
\qed

\begin{remark}
Let us recall that the construction in Remark 8.3 of \cite{pde}
 provides  a recipe that attaches  to any
$\D$-prolongation sequence $S^*$ of Noetherian $p$-adically complete
rings over $A$, and to any pair $(E,\omega)$ consisting of an
elliptic curve $E/S^0$ and an invertible $1$-form $\omega$ on $E$,
elements
$$f^s_p(E,\omega,S^*), f^s_q(E,\omega,S^*) \in S^s$$
for all $s \geq 1$.
(Again, we will usually drop $S^*$ from notation.) The
rules $f^s_p$ and $f^1_q$ coincide with the forms introduced in
Theorem \ref{gama},  so they are, in particular, $\D$-modular forms.
However, notice that for $s \geq 2$, the
rules $f^s_q$ are not  $\D$-modular forms. For instance, for
$s=2$, we have the transformation law
$$f^2_q(E,\lambda \omega)=\lambda^2
f^2_q(E,\omega)+2 \lambda (\dd \lambda)f^1_q(E,\omega)\, .$$
\qed
\end{remark}

\subsection{$\D$-Fourier expansions}

Assume, in the discussion below, that $N>4$.
 We start by recalling the background of {\it classical} Fourier expansions;
we will freely use the notation in  \cite{DI}, p. 112. (The
discussion in  \cite{DI}, p. 112, involves the model $X_{\mu}(N)$
instead of the model $X_1(N)$ used here, but these two models, and
hence the two theories, are isomorphic over $\ZN$ cf. \cite{DI}, p.
113.)

There is a point $s_{\infty}:\ZN \ra X_1(N)_{\ZN}$ arising
from the generalized elliptic curve ${\mathbb P}^1_{\ZN}$ with its
canonical embedding of $\mu_{N,\ZN} \simeq ({\mathbb Z}/N{\mathbb
Z})_{\ZN}$; the complex point corresponding to $s_{\infty}$ is the
cusp $\Gamma_1(N) \cdot \infty$. The map $s_{\infty}$ is a closed
immersion. We denote by $\tilde{X}_1(N)_{\ZN}$ the completion of
$X_1(N)_{\ZN}$ along the image of $s_{\infty}$.

Let us now consider an indeterminate $\bq$ and the Tate generalized elliptic
curve
$$
{\rm Tate}(\bq)/\ZN [[\bq]]\, ,
$$
defined as the projective closure of the plane affine curve given by
$$
y^2=x^3-\frac{1}{48}E_4(\bq)x-\frac{1}{864}E_6(\bq)\, ,
$$
where $E_4$ and $E_6$ are the Eisenstein series
$$\begin{array}{rcl}
E_4(\bq) & = & 1+240 \cdot s_3(\bq)\, ,\\
E_6(\bq) & = & 1-504 \cdot s_5(\bq)\, ,
\end{array}
$$
and $s_m$ is defined by
$$
s_m(\bq):=\sum_{n \geq 1} \frac{n^m \bq^n}{1-\bq^n} \in \bZ[[\bq]]\, .
$$
(The normalization chosen here for the Tate curve is the same as
that in \cite{difmod}, but slightly different from the usual one, for
instance, that used in  \cite{pde}; the Tate curve used in
the present paper is, however,  isomorphic over
$\bZ[\sqrt{-1}][[\bq]]$, and hence over $R[[\bq]]$, to the Tate
curve in \cite{pde}.)
 The curve ${\rm Tate}(\bq)$ is equipped
with a canonical immersion
$$ \alpha_{can}: \mu_{N,\ZN} \simeq ({\mathbb
Z}/N{\mathbb Z})_{ \ZN} \rightarrow {\rm Tate}(\bq)\, ,$$
so there is an induced map ${\rm Spec} \,
\ZN[[\bq]] \ra X_1(N)_{\ZN}$. There is an induced isomorphism
\begin{equation}
\label{ttt}
{\rm Spf}\; \ZN[[\bq]] \ra \tilde{X}_1(N)_{\ZN}\, .
\end{equation}
There is also a canonical $1$-form $\omega_{can}=dx/y$ on the elliptic
curve ${\rm Tate}(\bq)$ over
$$\ZN ((\bq)):=\ZN [[\bq]][1/\bq]\, ,$$ such that,
for any classical modular form $f$ on $\Gamma_1(N)$ over $\bC$, the
series
$$f_{\infty}:=f_{\infty}(\bq):=
f({\rm Tate}(\bq)/\bC((\bq)),\alpha_{can}, \omega_{can})$$
 has image in
$\bq\bC[[\bq]]$, and is the classical Fourier expansion at the cusp
$\Gamma_1(N) \cdot \infty$. (Here we view $f$ as a function of
triples in the sense of \cite{Katz}.)

 Let us now move onto the $\D$-theory. We fix a
prime $p$ that does not divides $N$, a $\D$-ring $A$ that is an integral
$p$-adically complete Noetherian domain of characteristic zero, and a
homomorphism $\ZN \ra A$. Let $\bq^{(i,j)}$ be  indeterminates over $A$
parameterized by non-negative integers $i$, $j$, such that
$\bq^{(0,0)}=\bq$. We set $A((\bq)):=A[[\bq]][\bq^{-1}]$, and
$$S^n_{pq,\infty}:=A((\bq))\wh[\bq^{(i,j)}\mid_{1 \leq i+j \leq n}]\wh \, .$$
 There is a unique structure of $\D$-prolongation sequence on
$S^n_{pq,\infty}$ over $A$ that extends that of
$A[\bq^{(\leq r)},\bq^{-1}]\wh$,
and sends $A[[\bq]]$ into $A[[\bq]][\d \bq, \dd \bq]\wh$.

 Finally, we  define the
{\it $\D$-Fourier expansion map}
$$
\begin{array}{rcl}
E_{\infty}:{M}^r_{pq}(\Gamma_1(N),A,w) & \ra & S_{pq,\infty}^r \\
f & \mapsto &
E_{\infty}(f)=f_{\infty}
\end{array}
$$
by the formula
$$
f_{\infty}:=f_{\infty}(\bq^{(i,j)}|_{0 \leq i+j\leq r}):=
f({\rm Tate}(\bq)/S_{pq,\infty}^0,\alpha_{can},\omega_{can},
S_{pq,\infty}^*)\, .
$$
By Lemma \ref{allfa}, we get an induced map
$$
E_{\infty}:{M}^r_{pq}(A,w) \ra S^r_{pq,\infty}\, .
$$
Since the normalized Eisenstein series $E_{p-1}(q)$ is congruent to
$1$ mod $p$, we get $\D$-Fourier expansion maps
\begin{equation}
\begin{array}{c}
E_{\infty}:{M}^r_{pq,ord}(\Gamma_1(N),A,w) \ra S^r_{pq,\infty}
\, , \\
E_{\infty}:{M}^r_{pq,ord}(A,w) \ra S^r_{pq,\infty}\, .
\end{array}
\label{fem}
\end{equation}

\begin{proposition}
\label{deltaexp} {\rm (}$\D$-expansion principle{\rm )} The $\D$-Fourier
expansion maps {\rm (\ref{fem})} are injective, with torsion free cokernel.
\end{proposition}

{\it Proof}.
This follows exactly as in  \cite{book}, Proposition 8.29, where the
Serre-Tate expansions (rather than the Fourier expansions) were
considered. \qed
\medskip

If instead of $\D$ we use $\d$ or $\dd$) only, we arrive at a
$\d$-prolongation sequence
$$S^n_{p,\infty}:=A((\bq))\wh[\bq^{(i,0)}: \, 1\leq i \leq n]\wh\, ,$$
or at a $\dd$-prolongation sequence
$$S^n_{q,\infty}:=A((\bq))[\bq^{(0,j)}: \, 1\leq j \leq n]\, {\rm )}\, ,$$
respectively.
We have corresponding $\d$-Fourier or $\dd$-Fourier expansion
maps. These ``ordinary'' objects are compatible with the partial
differential objects in an obvious way.

We introduce now certain subspaces of $S_{pq,\infty}^r$. For each
prime  $l$  different from $p$, there is a unique morphism of
$\D$-prolongation sequences $\varphi_l^*:S_{pq,\infty}^* \ra
S_{pq,\infty}^*$ such that $\varphi_l^0(F(\bq))=F(\bq^l)$. Then, for each
even integer $m$, we define the following $A$-submodules of
$S_{pq,\infty}^r$:
$$I^r_{pq,\infty}(A,m):=\{f \in S_{pq,\infty}^r\ |\
\varphi_l(f)=l^{-m/2} \cdot f\ \ \text{for all $l \neq p$}\}\, .$$

\begin{proposition}
\label{subb} Suppose that $w$ is a weight of even degree $m$.
Then the $\D$-Fourier expansion maps {\rm (\ref{fem})}
{\rm (}see Proposition {\rm \ref{deltaexp})} send both,
$I^r_{pq}(\Gamma_1(N),A,w)_1$ and $I^r_{pq,ord }(\Gamma_1(N),A,w)_1$,
into $I^r_{pq,\infty}(A, m)$.
\end{proposition}

{\it Proof}. We use the same argument as in the proof of
Proposition 7.13 in \cite{difmod}. \qed

\subsection{Convection equations}
In what follows, our purpose is to investigate the space
$I^1_{pq,ord}(A,-2)$.

Consider the following two series in $S^1_{pq,\infty}$:
$$\begin{array}{rcl}
 \Psi_p & = & {\displaystyle \frac{1}{p} \cdot \log{\left(
1+p\frac{\d \bq}{\bq^p} \right)}=\sum_{n=1}^{\infty} (-1)^{n-1}
p^{n-1} n^{-1} \left( \frac{\d \bq}{\bq^p} \right)^n }\, , \vspace{1mm} \\
\Psi_q & = & {\displaystyle \frac{\dd \bq}{\bq}}\, .
\end{array}
$$
Clearly we have $\Psi_p,\Psi_q \in I^1_{pq,\infty}(A,-2)$.

Let $L$ be the fraction field of $A$.

\begin{proposition}
\label{bass} The series $\Psi_p,\Psi_q$ form an $L$-basis of
$I^1_{pq,\infty}(A,-2) \otimes L$.
\end{proposition}

{\it Proof}. Consider the map
$$\begin{array}{rcl}
c:I^1_{pq,\infty}(A,-2) & \ra &  A^2 \\
c(f) & = & (c_p(f),c_q(f))
\end{array}\, ,
$$
where $c_p(f)$ and $c_q(f)$ are the
coefficients of $\frac{\d \bq}{\bq^p}$ and  $\frac{\dd
\bq}{\bq}$ in $f$, respectively. It is enough to show that the map $c$ in
injective.

For any $i,j \geq 0$, let us set
$$
\partial_{ij}=\frac{\partial^{i+j}}{(\partial \d \bq)^i (\partial
\dd \bq)^j}\, ,
$$
and let $g=g(\bq,\d \bq,\dd \bq) \in
I^1_{pq,\infty}(A,-2)$. Since $\dd (\bq^2)=2 \bq \dd \bq$ and $\d
(\bq^2)= 2\bq^p \d \bq+p(\d \bq)^2$, it follows that
$$
g(\bq^2,  2\bq^p \d \bq+p(\d \bq)^2,2\bq\dd \bq)=\varphi_2(g)=2
g(\bq,\d \bq, \dd \bq)\, .
$$
We may easily check by induction that
\begin{equation}
\label{indd}
(\partial_{ij} g)(\bq^2,  2\bq^p \d
\bq+p(\d \bq)^2,2\bq\dd \bq)2^{i+j} \bq^i (\bq^p+p\d \bq)^j +U_{ij}=
2 \cdot
(\partial_{ij} g)(\bq, \d\bq, \dd \bq)\, ,
\end{equation}
where
$U_{ij}$ is a linear combination with coefficients in the ring
$S_{pq,\infty}^1$ of series of the form
$$(\partial_{i' j'} g)(\bq^2,  2\bq^p \d \bq+p(\d
\bq)^2,2\bq\dd \bq)$$ with $i'+j'<i+j$. If we assume that $g$ is in the kernel
of the map $c$, setting $\dd \bq=\d \bq=0$ in (\ref{indd}), and using
Lemma 7.22 in \cite{difmod}, we then obtain by induction that
$(\partial_{ij}g)(\bq,0,0)=0$ for all $i,j$, and so $g$ must be identically
$0$.
\qed.
\medskip

In what follows, our  $\D$-ring $A$ will be the ring
$R=(\bZ_p^{ur})\wh$. We shall denote by $K=R[1/p]$ its fraction field.

\begin{lemma}
\label{cafgata}
There exist $c, \gamma \in \bZ_p^{\times}$ such that
the following hold:
\begin{enumerate}
\item $f^1_p$ has $\D$-Fourier expansion $c \cdot \Psi_p$.
\item $f^1_q$ has $\D$-Fourier expansion $\gamma \cdot \Psi_q$.
\item $f^{\partial}_p$  has $\D$-Fourier expansion
$1$.
\end{enumerate}
\end{lemma}

{\it Proof}.
 Cf. \cite{difmod}, Corollary
7.26, and \cite{barcau},  Construction 3.2 and Theorem 5.1. \qed
\medskip

Let us recall the form $f^1_{p,-2}=f^{\partial}_pf^1_p$
in Remark \ref{cano}.

\begin{theorem}
\label{main} The forms $f^1_q$ and $f^1_{p,-2}$ form a
$K$-basis of $I^1_{pq,ord}(R,-2)\otimes K$.
\end{theorem}

In particular, any element $f \in I^1_{pq,ord}(R,-2)$ can be written
as  a $K$-linear combination
\begin{equation}
\label{ace}
f=\varphi^1_q + \varphi^1_p\, ,
\end{equation}
where $\varphi^1_q \in K \cdot f^1_q$ and $\varphi^1_p \in K \cdot f^1_{p,-2}$.
Such a linear combination can be referred to as an
{\it arithmetic convection equation}.

{\it Proof.} By Propositions \ref{deltaexp} and \ref{subb}, we have
an injective map
$$I^1_{pq,ord}(\Gamma_1(N),R,-2)_1\otimes K \ra
I^1_{pq,\infty}(R,-2) \otimes K.$$ The result follows by
Lemma \ref{cafgata} and Proposition \ref{bass}.
\qed
\medskip

If we combine the Theorem above and the observations in Remark \ref{cano},
we obtain bases for any of the spaces $I^1_{pq,ord}(R,w)$ with $ord(w)\leq 1$
and $deg(w)=-2$.

We investigate the solution spaces of equations defined by this type of
forms in the next section.

\subsection{Heat equations}
Motivated by (\ref{ace}) and the notion introduce above, we now
define an {\it arithmetic heat equation} to be one given by
a form of the type
\begin{equation}
\label{ahe}
f=\varphi^1_q+\varphi^2_p \in I^2_{pq,ord}(R,-2)\, ,
\end{equation}
where $\varphi^1_q \in I^1_q(R,-2)$
and $\varphi^2_p \in I^2_{p,ord}(R,-2)$. Thus, arithmetic convection
equations are special cases of arithmetic heat equations, those where
$\varphi^2_p \in I^1_{p,ord}(R,-2)$.

We investigate their solution sets in the next section also.

\subsection{Non-existence of wave equations}
Motivated once again by (\ref{ace}) and now (\ref{ahe}) also, we could try
to introduce arithmetic analogues of wave equations by looking at
forms of the type
$$
f=\varphi^2_q+ \varphi^2_p \in I^2_{pq,ord}(R,-2)\, ,
$$
where $\varphi^2_q \in I^2_q(R,-2)$ and $\varphi^2_p \in
I^2_{p,ord}(R,-2)$, and the corresponding equations they define. However,
our next result shows that doing so does not lead to anything genuinely new.

\begin{theorem}
\label{incauna} We have that $I^1_q(R,-2)=I^2_q(R,-2)$.
\end{theorem}

\begin{remark}
We  expect that $I^1_q(R,-2)=I^r_q(R,-2)$ for all $r \geq 2$. However,
the arguments that we will use to prove Theorem \ref{incauna}
seem to work for the case $r=2$ only. \qed
\end{remark}

\begin{remark}
We expect that $I^2_{pq,ord}(R,-2)\otimes K$ has a basis consisting
of the elements
\begin{equation}
\label{inimioara}
f^1_q,\ (f^{\partial}_p)^2(f^1_q)^{\phi},\
f^{\partial}_pf^1_p,\ (f^{\partial}_p)^{\phi+2}
(f^1_p)^{\phi}\, ,
\end{equation}
which would represent a strengthening of Theorem \ref{incauna}.
Indeed, if this Theorem were
false, then $I^2_q(R,-2) \otimes K$ would have dimension $\geq 2$,
and so two $K$-linearly independent elements of this space together
with the last $3$ of the $4$ elements in (\ref{inimioara})
would yield $5$ linearly independent elements in
$I^2_{pq,ord}(R,-2)\otimes K$. \qed
\end{remark}

\medskip
In order to prove Theorem \ref{incauna}, we now discuss some preliminary
results. The strategy that we follow in this proof is inspired by a
method of Barcau \cite{barcau}. For convenience, we shall use the notation
$x^{(r)}$ for $\dd^r x$, and when appropriate, the notation
$x^{'}$ and $x^{''}$ for $x^{(1)}$ and $x^{(2)}$, respectively.

First, we recall \cite{difmod} that the space
$M^r_q(R,m)$, $m \in \bZ$,  identifies with the space of all
rational functions
$$
f \in M^r_q:=R[a_4,a_6,\ldots ,a_4^{(r)},a_6^{(r)},\Delta^{-1}]
$$
that have {\it weight} $m$, in the sense that
\begin{equation}
\label{gain}
f(\ldots ,\dd^i (\Lambda^4 a_4), \dd^i(\Lambda^6
a_6),\ldots ,\Lambda^{12}\Delta^{-1})=\Lambda^m f(\ldots ,a_4^{(i)},
a_6^{(i)},\ldots ,\Delta^{-1})
\end{equation}
in the ring
$$
R[a_4,a_6,\Lambda,\ldots ,a_4^{(r)},a_6^{(r)},\Lambda^{(r)}, \Lambda^{-1},
\Delta^{-1}]\, .
$$
We recall also \cite{difmod} that we have natural maps ---called the
$(\bq,\ldots ,\bq^{(r)})$-Fourier expansion maps---
\begin{equation}
\label{room}
\begin{array}{rcl}
M^r_q & \rightarrow & R((\bq))[\bq',\ldots ,\bq^{(r)}]\\
u & \mapsto & u_{\infty}
\end{array}\, .
\end{equation}

We have the following:

\begin{lemma}
\label{x5} The $(\bq,\bq')$-Fourier expansion map $M^1_q \ra
R((\bq))[\bq']$ is injective.
\end{lemma}

{\it Proof.} This is contained in \cite{difmod}, Proposition 7.10.
\qed

\begin{remark}
Notice that the map in (\ref{room}) is not injective for $r
\geq 2$; cf. \cite{difmod}, Proposition 7.10. \qed
\end{remark}

\medskip

On the ring $M^r_q$, we introduce the following derivation operators:
$$
\begin{array}{rcl}
\partial_r & = & 16a_4^2\frac{\partial}{\partial a_6^{(r)}}-72 a_6
\frac{\partial}{\partial a_4^{(r)}}\, , \vspace{1mm} \\
D_r & = & 4a_4\frac{\partial}{\partial a_4^{(r)}}+
6a_6\frac{\partial}{\partial a_6^{(r)}}\, .
\end{array}
$$
Notice that $\partial_0$ is the usual Serre operator (see \cite{Katz},
Appendix), while $D_0$ is the usual Euler operator for weighted
homogeneous polynomials. In general, these operators
are ``geometric'' analogues of certain ``arithmetic'' operators
introduced in \cite{barcau}.

\begin{lemma}
\label{x1}
If $f \in M^r_q$ has weight $m$, then $\partial_r f$ has
weight $m+2$.
\end{lemma}

{\it Proof}. Apply the operator $\partial_r$ in (\ref{gain}). \qed

\begin{lemma}
\label{x2} If $f \in M^r_q$ has weight $m$ and $r \geq 1$, then $D_r
f=0$.
\end{lemma}

{\it Proof}. Apply the operator $\frac{\partial}{\partial \Lambda^{(r)}}$ in
(\ref{gain}), and set $\Lambda=1$. \qed

\begin{lemma}
\label{x3} If $f \in M^r_q$ has weight $m$ and $r \geq 1$, then we have
the equality of $\dd$-Fourier series
$$
(\partial_r f)_{\infty}=12 \bq \frac{\partial f_{\infty}}{\partial
\bq^{(r)}}\, .
$$
\end{lemma}

{\it Proof}. Let us write $(a_4)_{\infty}=a_4(\bq)$, with a similar
notation for $a_6$. We have that
\begin{equation}
\begin{array}{rcl}
12 \bq \frac{d a_4(\bq)}{d \bq} & = 4P(\bq) a_4(\bq) -72
a_6(\bq)\, ,\vspace{1mm} \\
12 \bq \frac{d a_6(\bq)}{d \bq} & = 6P(\bq) a_6(\bq) +16 a_4^2(\bq)\, ,
\end{array}
\end{equation}
where $P(\bq)\in R[[\bq]]$ is the Ramanujan series (\cite{Katz},
Appendix). Consequently, if $V$ denotes the vector
$$
(\ldots ,\dd^i(a_4(\bq)),\dd^i(a_6(\bq)),\ldots )\, ,
$$
we have that
$$
\begin{array}{rcl}
12 \bq \frac{\partial f_{\infty}}{\partial \bq^{(r)}} & = & 12\bq
\left[ \frac{\partial f}{\partial
a_4^{(r)}}(V)\frac{\partial}{\partial
\bq^{(r)}}(\dd^r(a_4(\bq)))+\frac{\partial f}{\partial
a_6^{(r)}}(V)\frac{\partial}{\partial \bq^{(r)}}(\dd^r(a_6(\bq)))\right]
\vspace{1mm} \\
 & = & 12\bq \left[ \frac{\partial f}{\partial
a_4^{(r)}}(V)\frac{d a_4(\bq)}{d \bq}+\frac{\partial f}{\partial
a_6^{(r)}}(V)\frac{d
a_6(\bq)}{d \bq}\right] \vspace{1mm} \\
 & = & (4P(\bq) a_4(\bq)\! - \! 72 a_6(\bq))\frac{\partial f}{\partial
a_4^{(r)}}(V)\! + \! (6P(\bq) a_6(\bq) +16 a_4^2(\bq))\frac{\partial
f}{\partial a_6^{(r)}}(V) \vspace{1mm} \\
 & = & P(\bq) (D_r f)_{\infty}+ (\partial_r f)_{\infty}\, .
\end{array}
$$
The desired result now follows by Lemma \ref{x2}.
\qed

\begin{lemma}
\label{x4}
We have that $\dd f^1_q \not\in M^2_q(R,-2)$.
\end{lemma}

{\it Proof}. Let us assume that the opposite is true, that is to say,
that
$$
(\dd f^1_q)(E,\lambda\omega,S^*)=\lambda^{-2}(\dd f^1_q)(E,\omega,S^*)\, .
$$
We apply $\dd$ in the identity
$$
f^1_q(E,\lambda\omega,S^*)=\lambda^{-2}f^1_q(E,\omega,S^*)\, ,
$$
and use the assumption to get that
$$
\lambda (\dd \lambda) f^1_q(E,\omega,S^*) =0\, ,$$
which is, of course, false for the generic $(E,\omega)$ and
for any $\lambda$ with $\dd \lambda \neq 0$.
\qed
\medskip

We will also need the following explicit formula for $f^1_q$ that is due to
Hurlburt \cite{hurl}; cf also  \cite{difmod}, Proposition 4.10 and Corollary
7.24.

\begin{lemma}
\label{x6} \cite{hurl}
We have $f^1_q=\gamma_1 \frac{2a_4a_6'-3a_6a'_4}{\Delta}$, where $\gamma_1
\in R^{\times}$.
\end{lemma}

(Of course the constant $\gamma_1$ above and the constant $\gamma$
in Lemma \ref{cafgata} are explicitly related to each other, but this relation
is irrelevant to us here.)

 {\it Proof of Theorem {\rm \ref{incauna}}}.
Let $g \in I^2_q(R,-2)$. We want to show that $g \in I^1_q(R,-2)$.
By \cite{difmod}, Proposition 7.13 and 7.23, we may assume that
\begin{equation}
\label{tz}
g_{\infty}=\lambda \frac{\bq'}{\bq}+ \mu \left(
\frac{\bq}{\bq}\right)' \, ,
\end{equation}
with $\lambda, \mu \in R$. Replacing $g$ by $g$ plus a multiple of
$f^1_q$, we may assume that $\lambda=0$.
We will prove that $g=0$.

By Lemma \ref{x3}, we have that
$$
(\partial_2 g)_{\infty}=12\bq \frac{\partial g_{\infty}}{\partial
\bq''}=12 \mu \, .
$$
By Lemma \ref{x1}, we see that $\partial_2g$ has weight $0$. Now,
the constant $12 \mu \in M^2_q$ has weight $0$ also, and
$\dd$-Fourier expansion $(12\mu)_{\infty}=12 \mu$. By Proposition
\ref{deltaexp}, it follows that
\begin{equation}
\label{y1}
\partial_2 g=12 \mu\, .
\end{equation}
Notice that by Lemma \ref{x2} we have that
\begin{equation}
\label{y2} D_2 g=0\, .
\end{equation}
By Lemma \ref{x6}, we have
$$\dd f^1_q \in \gamma_1 \frac{2a_4a_6''-3a_6a''_4}{\Delta}+M^1_q\, .$$
A trivial computation yields
\begin{equation}
\label{z1}
\partial_2 \dd f^1_q=\-\gamma_1/2\, , \quad
D_2 \dd f^1_q=0 \, .
\end{equation}
We set $h:=g+\frac{24 \mu}{\gamma_1} \dd f^1_q \in M^2_q$. By
(\ref{y1}), (\ref{y2}) and (\ref{z1}), we obtain that
\begin{equation}
\label{w1}
\partial_2 h=0\, ,\quad D_2h=0\, .
\end{equation}
Also, by (\ref{tz}), $h$ has $(\bq,\bq')$-Fourier expansion
\begin{equation}
\label{a1}
h_{\infty}=\gamma_2 \left( \frac{\bq'}{\bq} \right)'\, ,\quad
\gamma_2=\mu+\frac{24\mu \gamma}{\gamma_1}\, .
\end{equation}
Now, by (\ref{w1}),  we get that
$$
\frac{\partial h}{\partial a_4''}=\frac{\partial h}{\partial
a_6''}=0\, ,
$$
hence $h \in M^1_q$. On the other hand, the element
$\frac{\gamma_2}{\gamma} \dd f^1_q \in M^1_q$ has
$(\bq,\bq')$-Fourier expansion
\begin{equation}
\label{a2} \left(\frac{\gamma_2}{\gamma} \dd f^1_q \right)_{\infty}=
\gamma_2 \left( \frac{\bq'}{\bq} \right)'\, .
\end{equation}
By Lemma \ref{x5} and (\ref{a1}) and (\ref{a2}), we must have that
$h=\frac{\gamma_2}{\gamma}\dd f^1_q$, and therefore,
$g=\gamma_3 \dd f^1_q$, $\gamma_3 \in R$. By Lemma \ref{x4}, we
conclude that $\gamma_3=0$, hence $g=0$, and we are done. \qed

\section{Solutions}
\setcounter{theorem}{0}
We now analyse the sets of solutions
of the equations introduced in the previous section.
We first examine solutions with bad reduction at $q=0$; then we shall
examine solutions with good reduction at $q=0$.

\subsection{Bad reduction}
We define the following subset of $\bM(R((q))\wh)$:
$$\bM(R((q))\wh)_{bad}:=\{(a,b) \in R[[q]]^{\times} \times
R[[q]]^{\times};\ 4a^3+27b^2 \in qR[[q]]^{\times}\} \subset
\bM(R((q))\wh)\, .$$
We will prove the following

\begin{lemma}
\label{capitula} We have that
$\bM(R((q))\wh))_{bad} \subset \bM_{ord}(R((q))\wh)$.
\end{lemma}

We assume this Lemma for the time being. For any $f
\in M^r_{pq,ord}(R,-2)$, we define the set of {\it bad reduction solutions}
of the equation $f=0$ by
$$\cU_{bad}:=\{(a,b) \in \bM(R((q))\wh)_{bad}\ ;\ f(a,b)=0\}\, .$$
We aim at the description of this space for arithmetic heat equations
$$
f:=f^1_q+\varphi^2_p\, ;
$$
cf. (\ref{ahe}).

By Theorem \ref{gama},
$$
\varphi^2_p=\lambda  f^1_{p,-2}+\epsilon f^{1 \phi}_{p,-2}\, ,
$$
where $\lambda,\epsilon \in K$.
In studying the equation $f=0$, we may assume that
$\lambda,\epsilon \in R$. We will further assume that
$\lambda \in R^{\times}$.
Let
$$
\int:K[[q]]\ra qK[[q]]
$$
be the usual integration operator. Also consider the Fourier expansions
of $a_4,a_6$ respectively:
$$
a_{4,\infty}(\bq):=-\frac{1}{48}E_4(\bq)\, , \quad
a_{6,\infty}(\bq):=-\frac{1}{864}E_6(\bq)\in R[[\bq]]\, .
$$
For $n \geq -1$, we define inductively the rational functions
$$b_n=b_n(x,y) \in \bZ[x,y,\frac{1}{(1-y)(1-y^2)\cdots (1-y^n)}]$$
by setting $b_{-1}=0$, $b_0=1$, and the recurrence relation
\begin{equation}
\label{recurrence}
b_n =\frac{1-x}{1-y^n} b_{n-1}+\frac{x}{1-y^n} b_{n-2}\, .
\end{equation}
Notice that we have
$$
\begin{array}{rcl}
 b_n(x,0) & = & {\displaystyle
\frac{1+(-1)^nx^{n+1}}{1+x}=1-x+x^2-\cdots +(-1)^nx^n}\vspace{1mm} \, ,\\
 b_n(0,y) & = & {\displaystyle \frac{1}{(1-y)(1-y^2)\cdots (1-y^n)}}\, .
\end{array}
$$
Both formulae above have, of course, a quantum theoretic flavor.

For $0 \neq \kappa \in \bZ_+$, $z \in \bQ_p$, $\alpha \in K$,
$\eta \in K^{\times}$ and $v \in K[[q]]^{\times}$, let us define
\begin{equation}
\label{ursu}
\begin{array}{rcll}
u_{a,\kappa,\alpha}^{z} & := & \sum_{n \geq 0} b_n(pz,p)
\alpha^{\phi^{n}}  q^{\kappa p^{n}} & \in K[[q]]\, , \vspace{1mm} \\
u_{m,\kappa,\alpha}^{z} & := &{\displaystyle  \exp{\left( \int
u^{z}_{a,\kappa,\alpha}\frac{dq}{q} \right)}} & \in 1+qK[[q]]\, ,
\vspace{1mm} \\
 u_{\Gamma,\eta, v}^{z} & := & (v^4a_{4,\infty}(\eta
q),v^6a_{6,\infty}(\eta q)) & \in \bM(K((q)))\, , \vspace{1mm} \\
u^{z}_{\Gamma,\eta,v,\kappa,\alpha} & := &  (v^4a_{4,\infty}(\eta q
u^{z}_{m,\kappa,\alpha}), v^6a_{6,\infty}(\eta q u^{z}_{m,\kappa,\alpha}))&
\in \bM(K((q)))\, .
\end{array}
\end{equation}
Here the indices $a,m,\Gamma$ stand for the additive group
${\mathbb G}_a$, the multiplicative group ${\mathbb G}_m$, and the
modular scheme $\bM_{\Gamma}$ attached to $\Gamma={\rm SL}_2(\bZ)$.
 As observed in \cite{pde}, the series
$$
u^{0}_{m, \kappa,\alpha}=\exp \left(\frac{\alpha q^{\kappa}}{\kappa}+
\sum_{n \geq 1}  \frac{\alpha^{\phi^n} q^{\kappa p^n}}{\kappa
p^n (1-p)(1-p^2)\cdots (1-p^n)}
\right)$$
 is a sort of hybrid between the Artin-Hasse exponential
(\cite{koch}, p. 138) and a quantum exponential (\cite{kac}, p. 30).

\begin{lemma}
\label{integraliti}
If $\kappa \not\in p\bZ$, $z \in \bZ_p$, and  $\alpha \in R$, then
$$
u_{a,\kappa,\alpha}^{z} \in R[[q]]\, , \quad
u_{m,\kappa,\alpha}^{z} \in 1+qR[[q]]\, .
$$
If, in addition, $\eta \in R^{\times}$ and $v \in R[[q]]^{\times}$, then
$$
u_{\Gamma,\eta, v}^{z}, u_{\Gamma,\eta, v,\kappa,\alpha}^{z}
\in \bM((R((q)))\, .
$$
\end{lemma}

{\it Proof.} We freely use the theory in \cite{pde}. Let us consider
the $\D$-character of $\bG_m$ defined by
$$
\psi_m:=\psi_1+(\kappa+\kappa z\phi)\psi_p\, ,
$$
where $\psi_q$, $\psi_p$ are defined by the series $\Psi_q$, $\Psi_p$,
respectively;
cf. \cite{pde}, Definition 7.3 and Example 7.4.
Using the terminology in loc. cit., we have that $\psi_m$
is non-degenerate, and its characteristic polynomial is
$$
\mu(\xi_p,\xi_q)=\xi_q+p \kappa z \xi_p^2+(\kappa-p \kappa z) \xi_p-\kappa\, .
$$
Thus, $\mu$ has characteristic integer $\kappa$.
By \cite{pde}, Definitions 7.3 and 6.2,
the basic series of $\psi_m$ are equal to the series $u^z_{m,\kappa,\alpha}$.
By \cite{pde}, Lemma 7.6, we have that
$u^z_{\kappa,\alpha} \in 1+qR[[q]]$. The remaining
portion of the Lemma is simple.
\qed

\begin{remark}
\label{integralitii}
Still using the notation above, we recall from \cite{pde}, Theorem 1.10,
that the set of solutions in $1+qR[[q]]$ of $\psi_m=0$ equals
$\{u^z_{m,\kappa,\alpha}\ |\ \alpha \in R\}$. \qed
\end{remark}

\begin{theorem}\label{main2}
Consider the arithmetic heat equation given by
$$
f:=f^1_q+\lambda  f^1_{p,-2}+\epsilon f^{1 \phi}_{p,-2}\, ,
$$
where
$\lambda \in R^{\times}$, $\epsilon \in pR$, and
$z:=\frac{\epsilon}{\lambda} \in p\bZ_p$.
Let $\kappa :=\frac{c}{\gamma} \cdot \lambda\in R^{\times}$, and
$\beta \in R^{\times}$ be such that
$$
\Psi_p(\beta)+z \Psi_p(\beta)^{\phi}=-\frac{1}{\kappa}\, .
$$
{\rm (}such a $\beta$ always exists{\rm )}.
If $\cU_{bad}$ is the set of bad reduction solutions of the equation
$f=0$, then the following hold:
\begin{enumerate}
\item Assume $\kappa \not\in {\mathbb Z}_+$. Then
$$
\cU_{bad}=\{u^z_{\Gamma,\zeta \beta, v}\, |\; \zeta \in \mu(R),\; v \in
R[[q]]^{\times}\}\, .
$$
\item Assume $\kappa \in {\mathbb Z}_+$. Then
$$
\cU_{bad}=\{u^z_{\Gamma,\zeta \beta, v, \kappa, \alpha} \, |\;
\zeta \in \mu(R),\, v \in R[[q]]^{\times},\ \alpha \in R\}\, .
$$
\end{enumerate}
\end{theorem}

We may view this as a ``quantization'' result: $\cU_{bad}$ has jumps
exactly at the integral positive values of $\kappa$.

Morally speaking, this Theorem is about arithmetic
heat equations that are close to, or coincide with, a
given arithmetic convection equation. The case of an arithmetic
convection equation corresponds to $\epsilon=0$; the case
$0 \neq \epsilon \in pR$
corresponds to arithmetic heat equations different that yield in the
limit the arithmetic convection equation.
The problem of how solutions behave when $\epsilon \ra 0$ will be addressed
in the next section. It is worth mentioning that heat equations of the
form $f^1_q+ \lambda f^2_{p,-2} =0$, say, are not covered
by our Theorem; such equations are not ``close to arithmetic
convection equations.''

We derive a preliminary result before the proof. We will use it
to prove Lemma \ref{capitula}, which in turn, shall be used in the proof
of the Theorem.

Let us introduce the map
\begin{equation}
\label{lazi}
\begin{array}{rcl}
\iota:R[[q]]^{\times} \times (R[[q]]^{\times}/\{\pm 1\})  & \ra
& \bM(R((q))\wh) \\
(u,v) & \mapsto & (v^4a_{4,\infty}(uq),v^6a_{6,\infty}(uq))\, .
\end{array}
\end{equation}

\begin{lemma}
\label{thcan} The map
$\iota$ in {\rm (\ref{lazi})} is injective, and its image is
$\bM(R((q))\wh)_{bad}$.
\end{lemma}

{\it Proof}. Let us assume that $\iota(u_1,v_1)=\iota(u_2,v_2)$. Then
$j_{\infty}(u_1 q)=j_{\infty}(u_2q)$, where
$$
j_{\infty}(\bq)=1/\bq+744+\cdots
$$
is the Fourier expansion of the $j$-invariant. It follows that
$\frac{1}{j_{\infty}}(u_1q)=\frac{1}{j_{\infty}}(u_2q)$. Since
$\frac{1}{j_{\infty}(\bq)}=\bq+\cdots $ has a compositional inverse,
it follows that $u_1q=u_2q$, hence $u_1=u_2$ and $v_1^2=v_2^2$.
This proves the injectivity of $\iota$.

Let us now verify the identity of sets
${\rm Im}\, iota=\bM(R((q))\wh)_{bad}$. That
the first set is included in the second is clear. Conversely, let
$(a,b) \in \bM(R((q))\wh)_{bad}$. Then
$$
\frac{4a^3+27b^2}{a^3} \in qR[[q]]^{\times}\, .
$$
We let $\sigma \in qR[[q]]^{\times}$ be the compositional inverse of
$\frac{1}{j_{\infty}(\bq)} \in \bq R[[\bq]]^{\times}$,
nd set
$$
u=\frac{1}{q} \sigma\left(2^83^3\frac{4a^3+27b^2}{a^3} \right) \in
R[[q]]^{\times}\, .
$$
Then the elliptic curves $y^2=x^3+ax+b$ and
$y^2=x^3+a_{4,\infty}(uq)x+a_{6,\infty}(uq)$ have the same
$j$-invariant, and therefore, there exists an element $v$ in the algebraic
closure of the fraction field of $R((q))$ such that
$a=v^4a_{4,\infty}(uq)$ and $b=v^6a_{6,\infty}(uq)$, respectively.
We get that $v^2 \in R[[q]]^{\times}$, and consequently, by Hensel,  that
$v \in R[[q]]^{\times}$.
So $(a,b) \in {\rm Im} \, \iota$, which proves the desired inclusion. \qed
\medskip

{\it Proof of Lemma {\rm \ref{capitula}}}. By Lemma \ref{thcan}, any
element in $\bM(R((q))\wh)_{bad}$ is of the form $\epsilon(u,v)$.
Let $E_{p-1} \in {\mathbb Z}_p[a_4,a_6]$ be the polynomial whose
Fourier expansion is the Eisenstein series $E_{p-1}(q)$. Then
$$
E_{p-1}(a_{4,\infty}(uq),a_{6,\infty}
(uq))=E_{p-1}(uq) \equiv 1\; {\rm mod}\; p\, ,
$$
so $\epsilon(u,v)$ corresponds to an ordinary elliptic curve.
\qed
\medskip

{\it Proof of Theorem {\rm \ref{main2}}}.
 If we
compose a map $g:\bM_{ord}(R((q))\wh) \ra R((q))\wh$ with the map
$\iota$ in (\ref{lazi}), we obtain a map
$$
g_{\iota}=g \circ
\iota:R[[q]]^{\times} \times (R[[q]]^{\times}/\{\pm 1\}) \ra
R((q))\wh \, .
$$
By Lemma \ref{thcan}, $\iota$ induces a bijection between $\cU_{bad}$
and the set
$$
f_{\iota}^{-1}(0)
\subset R[[q]]^{\times} \times (R[[q]]^{\times}/\{\pm 1\})\, .
$$
For any $s \in R((q))\wh$, and abusing notation a bit, we write
$\Psi_p(s)$ in place of $\Psi_p(s,\d s)$, and similarly for $\Psi_q$. Notice
that
$$
\Psi_p(s_1 s_2)=\Psi_p(s_1) + \Psi_p(s_2)\, ,
$$
with a similar expression for $\Psi_q$ instead.
By Lemma \ref{cafgata}, we have that
$$
\begin{array}{l}
(f^1_p)_{\iota}(u,v) = c \cdot v^{-1-\phi} \cdot \Psi_p(uq)
=c \cdot v^{-1-\phi} \cdot (\Psi_p(u)\! +\! \Psi_p(q))=c \cdot v^{-1-\phi}
\cdot \Psi_p(u)\, , \vspace{1mm} \\
(f^1_{q})_{\iota}(u,v)  = \gamma \cdot v^{-2}
\cdot \Psi_q(uq) =\gamma \cdot v^{-2} \cdot
(\Psi_q(u)+\Psi_q(q))=\gamma \cdot v^{-2} \cdot
(\Psi_q(u)+1)\, , \vspace{1mm} \\
(f^{\partial}_p)_{\iota}(u,v) = v^{\phi -1},
\end{array}
$$
so we have that
$$
f_{\iota}(u,v)=0
$$
if, and only if,
$$
\psi_m(u):= \Psi_q(u)+ \kappa \Psi_p(u)+\kappa z \phi \Psi_p(u)=-1\, .
$$
The latter has the obvious solution $u=\beta$, and any other solution
in $R[[q]]^{\times}$
is of the form
$u=\beta \cdot u_+$, where $u_+\in R[[q]]$ satisfies the homogeneous equation
\begin{equation}
\label{nuvine}
\Psi_m(u_+)=0\, .
\end{equation}
Clearly $u_0:=u_+(0)$ satisfies the equation
\begin{equation}
\label{ha}
\Psi_p(u_0)+z \phi \Psi_p(u_0)=0\, .
\end{equation}
This is equivalent to the condition $\Psi_p(u_0)=0$
(because
if $\Psi_p(u_0) \neq 0$ then the terms in the left hand side of (\ref{ha})
have distinct valuations). So $u_0=\zeta \in \mu(R)$ is a root of unity,
and $u_+=\zeta u_1$, $u_1 \in 1+qR[[q]]$, $\Psi_m(u_1)=0$.
By the proof of Lemma \ref{integraliti} and Remark \ref{integralitii},
if $\kappa \in \bZ$ then
we must have $u_1=u^z_{m,\kappa,\alpha}$ for some $\alpha \in R$.
On the other hand, if $\kappa \not\in \bZ_+$ then, by \cite{pde},
Theorem 7.10, $u_1=1$. This finishes the proof.
\qed

\begin{remark}
Notice that, by Lemma \ref{thcan}, the maps
$$\begin{array}{rcl}
(\zeta,v) & \mapsto & u_{\Gamma,\zeta \beta, v}\\
(\zeta, v,\alpha) & \mapsto & u_{\Gamma,\zeta \beta, v, \kappa,
\alpha}
\end{array}$$
are injective. Notice also that we may write
\begin{equation}
\label{capra} \cU_{bad}=\bigcup_{\zeta \in \mu(R)} \cU_{\zeta}\, ,
\end{equation}
where
$$\cU_{\zeta} = \{u_{\Gamma,\zeta \beta, v} |\, v \in
R[[q]]^{\times}\}$$
or
$$\cU_{\zeta} = \{u_{\Gamma,\zeta \beta, v, \kappa, \alpha} |\, v
\in R[[q]]^{\times}, \alpha \in R\}\, ,
$$
if $\kappa \not\in \bZ_+$ or $\kappa \in \bZ_+$, respectively.
According to these two cases, we have corresponding canonical identifications
\begin{equation}
\label{ggica}
\cU_{\zeta}/R[[q]]^{\times} \simeq 0\quad \text{or}\quad
\cU_{\zeta}/R[[q]]^{\times} \simeq R\, ,
\end{equation}
where the latter of the two is gotten via the $\alpha$ parameter.
This picture will have an analogue in the case of good reduction, where
$\mu(R) \simeq k^{\times}$
will be replaced by $k\backslash \{0,1728,ss\}$, $k=R/pR$
viewed as the $j$-line parameterizing elliptic curves over $k$, and
$ss$ the supersingular points. \qed
\end{remark}

\subsection{Good reduction}
The analysis of the set of bad solutions
has an analogue in the case of good reduction, which we now explore.

Let $(A,B) \in \bM(R)$. We
say that $(A,B)$ is a $CL$-{\it point} (a canonical lift point) if
the elliptic curve $y^2=x^3+Ax+B$ has ordinary reduction
mod $p$, and is isomorphic over $R$ to the canonical lift of its
reduction mod $p$. We denote by $\bM(R)_{CL} \subset \bM_{ord}(R)$
the set of CL points. Also, let $f=f^1_q+\lambda
f^1_{p,-2}+\epsilon f^{1\phi}_{p,-2}$,
with $\lambda \in R^{\times}$, $\epsilon \in pR$,  and set
$$
\begin{array}{rcl}
\cU_{ord} & = & \{(a,b) \in \bM_{ord}(R[[q]])\ |\
f(a,b)=0\}\, ,\\
\cU_0 & = & \{(A,B) \in \bM_{ord}(R)\ |\
f(A,B)=0\} \subset \cU_{ord}\, .
\end{array}
$$
We refer to $\cU_0$ as the set of {\it stationary solutions} of $f$.

\begin{lemma}
\label{L1}
We have that:
\begin{enumerate}
\item $\cU_0=\bM(R)_{CL}$.
\item The map
$$\begin{array}{rcl}
\bM_{ord}(R[[q]]) & \ra & \bM_{ord}(R) \\
(a,b) & \mapsto & (A,B)=(a(0),b(0))
\end{array}
$$
induced by the map
$$
q \mapsto 0\, ,
$$
sends $\cU_{ord}$ onto $\cU_0$.
\end{enumerate}
\end{lemma}

{\it Proof}. Let us check (1). By Lemma \ref{x6}, $f^1_q$ vanishes
on $\bM(R)$. Now, a pair $(A,B) \in \bM(R)$ belongs to
$\bM(R)_{CL}$ if, and only if, $f^1_p(A,B)=0$; cf. \cite{book},
Propositions 7.18 and 8.89.
On the other hand, $f^1_p(A,B)=0$
if, and only if,
$$
\lambda f^{\partial}_p(A,B)f^1_{p}(A,B)+\epsilon f^{\partial}_p(A,B)^{\phi+2}
 f^1_{p}(A,B)^{\phi}=0
$$
(because if $f^1_p(A,B) \neq 0$, then the terms in the sum above
have distinct valuations).

For (2), let $(a,b) \in \cU_{ord}$, and $(A,B)=(a(0),b(0))$. We have
that
$$
0= f^1_q(a,b)_{|q=0}+\lambda f^{\partial}_p(a,b)_{|q=0}
f^1_{p}(a,b)_{|q=0}+\epsilon f^{\partial}_p(a,b)^{\phi+2}_{|q=0} f^1_{p}
(a,b)^{\phi}_{|q=0}\, .
$$
By Lemma \ref{x6},
$$
f^1_q=\gamma_1 \frac{2aq \frac{db}{dq}-3bq
\frac{da}{dq}}{\Delta(a,b)} \in qR[[q]]\, ,
$$
so $f^1_q(a,b)_{|q=0}=0$.
Also
$f^1_p(a,b)_{|q=0}=f^1_p(A,B)$, and
$f^{\partial}_p(a,b)_{|q=0}=f^{\partial}_p(A,B)$ by the
compatibility of $f^1_p$ and $f^{\partial}_p$ with the
$\d$-ring
homomorphism
$$\begin{array}{rcl}
R[[q]] & \ra & R \\
q & \mapsto & 0
\end{array}\, .
$$
Thus, $f(A,B)=0$, and the map
$q \mapsto 0$ induces a map $\cU_{ord}\ra \cU_0$ that is
the identity on $\cU_0$, therefore, surjective. \qed

\medskip

In the sequel, we extend  the Serre-Tate expansion of modular forms
\cite{Mori} and $\d$-modular forms
\cite{shimura}, \cite{serre} to $\D$-modular forms. Let $(A_0,B_0)
\in \bM_{ord}(k)$, $k=R/pR$, and $E_0$ be the (ordinary)
elliptic curve $y^2=x^3+A_0 x+B_0$. We fix
a $\bZ_p$-basis of the physical Tate module
$T_p(E_0)$. Then, using the formal universal deformation space
$R[[t]]$ of $E_0$, the universal elliptic curve $E_{def}$ over it,
and the natural $1$-form $\omega_{def}$ on it (cf. \cite{shimura},
pp 212-213), we can define a $\D$-{\it Serre-Tate expansion map}
$$
\begin{array}{rcl}
M^n_{pq} & \ra &  S^n_{pq,def}:=R[[t]][\d^i \dd^j t\ ;\ 1 \leq i+j
 \leq n]\wh \vspace{1mm} \\
f & \mapsto & f_{E_0}:=f(E_{def},\omega_{def},S^*_{pq,def})
\end{array}\, .
$$

\begin{remark}
\label{iepure} A change of basis in the Tate module $T_p(E_0)$
(corresponding to an element $\lambda \in \bZ_p^{\times}$) has the
effect of composing the map $f \mapsto f_{E_0}$ with the unique
automorphism of prolongation sequences of $S^n_{pq,def}$ that sends
$t$ into $[\lambda^2](t)$; the latter series is the multiplication
by $\lambda^2$ in the formal group of ${\bG}_m$. \qed
\end{remark}
\medskip

For the next result, we consider the following series:
$$
\begin{array}{rcl}
\Psi_q(1+t)& = & \frac{\dd t}{1+t} \in R[[t]][\dd t] \vspace{1mm} \\
\Psi_p(1+t) & = & \frac{1}{p} \log
\frac{\phi(1+t)}{(1+t)^p}=\frac{1}{p} \log \left(1+p
\frac{\d t +\frac{1+t^p-(1+t)^p}{p}}{(1+t)^p} \right)\in R[[t]][\d t]\wh\, .
\end{array}
$$

\begin{lemma}\label{L0}
There exist $c_{E_0},\gamma_{E_0} \in R^{\times}$ such that the
following hold:
\begin{enumerate}
\item $f^1_p$ has $\D$-Serre-Tate expansion $c_{E_0} \cdot \Psi_p(1+t)$.
\item $f^1_q$ has $\D$-Serre-Tate expansion $\gamma_{E_0} \cdot
 \Psi_q(1+t)$.
\item $f^{\partial}_p$ has $\D$-Serre-Tate expansion $1$.
\end{enumerate}
\end{lemma}

 {\it Proof}. (1) was proved in \cite{shimura}, Lemma 2.4. (2) follows
exactly as in \cite{shimura}, pp. 230-234. (3) was proved in \cite{serre},
Proposition 7.2.
\qed

\begin{remark}
It is not clear what the relation between the constants
$c_{E_0}$ and $\gamma_{E_0}$ in the Lemma above and the constants $c$
and $\gamma$ in Lemma \ref{cafgata} is. \qed
\end{remark}

\begin{remark}
\label{kopacab}
Let $a_{4,E_0}, a_{6,E_0} \in R[[t]]$ be the Serre-Tate expansions
of $a_4$ and $a_6$, respectively. Since the specialization
$$\begin{array}{rcl}
R[[t]] & \ra & R \\ t & \mapsto & 0
\end{array}
$$
sends $E_{def}$ into an elliptic curve $E/R$ that is isomorphic to the
canonical lift of $E_0/k$, we get that
$(a_{4,E_0}(0), a_{6,E_0}(0)) \in \bM(R)_{CL}$. \qed
\end{remark}
\medskip

Let us consider the sets
$$
\bM(R[[q]])_{CL,E_0} \subset \bM(R[[q]])_{E_0} \subset
\bM_{ord}(R[[q]])\, ,
$$
defined as follows:
\begin{enumerate}
\item $\bM(R[[q]])_{E_0}$ consists of all $(a,b) \in \bM(R[[q]])$
such that the reduction mod $p$ of the elliptic curve defined by
$y^2=x^3+a(0)x+b(0)$ is isomorphic over $k$ to $E_0$.
\item $\bM(R[[q]])_{CL,E_0}$ consists of all $(a,b) \in \bM(R[[q]])$
such that  the elliptic curve defined by
$y^2=x^3+a(0)x+b(0)$ is isomorphic over $R$ to the canonical lift of
$E_0$.
\end{enumerate}
Clearly
$$
\bM_{ord}(R[[q]])=\bigcup_{E_0/k \ \text{ordinary}}
\bM(R[[q]])_{E_0}\, .
$$

On the other hand, let us consider the map
\begin{equation}
\label{epy}
\begin{array}{rcl}
\iota_{E_0}:qR[[q]] \times
(R[[q]]^{\times}/\{\pm 1\}) & \ra & \bM(R[[q]])\\
(u,v) & \mapsto & (v^4a_{4,E_0}(u),v^6 a_{6,E_0}(u))
\end{array} \, .
\end{equation}

\begin{lemma}
\label{L2}
Assume $j(E_0) \neq 0,1728$ in $k$ {\rm (}that is to say, $A_0B_0 \neq
0${\rm )}. Then the map $\iota_{E_0}$ in {\rm (\ref{epy})} is
injective, and its image is equal to $\bM(R[[q]])_{CL,E_0}$.
\end{lemma}

{\it Proof}.  We prove that ${\rm Im}\, \iota_{E_0} =\bM(R[[q]])_{CL,E_0}$.
The inclusion of the first set into the second follows by Remark
\ref{kopacab}. Let us check the oposite inclusion. Assume
$(a,b) \in \bM(R[[q]])_{CL,E_0}$, and set
$A=a(0)$, $B=b(0)$. So $y^2=x^3+Ax+B$ is isomorphic over $R$ to
$y^2=x^3+a_{4,E_0}(0)x+a_{6,E_0}(0)$. Thus, $A=\lambda^4
a_{4,E_0}(0)$, $B=\lambda^6 a_{6,E_0}(0)$ for some $\lambda \in
R^{\times}$, and without loosing generality, we assume that $\lambda=1$.
Let $j_{E_0} \in R[[t]]$ be
the Serre-Tate expansion of the $j$-invariant $j \in
\bZ_p[a_4,a_6,\Delta^{-1}]$. Since the $j$-invariant is unramified
on the fine moduli spaces $Y_1(N)$ ($N \geq 4$) above points
$\not\equiv 0,1728$ mod $p$, it follows that
$\frac{dj_{E_0}}{dt}(0) \in R^{\times}$. So the series
 $j_{E_0}-j_{E_0}(0) \in tR[[t]]$ has a compositional inverse
 $\sigma(t) \in tR[[t]]$. Notice that $j_{E_0}(0)=j(A,B)$. Let
$u:=\sigma(j(a,b)-j(A,B))\in qR[[q]]$. Then
$$
j_{E_0}(u)-j_{E_0}(0)=j(a,b)-j(A,B)\, ,
$$
hence $j_{E_0}(u)=j(a,b)$. Thus, the elliptic curves $y^2=x^3+ax+b$
 and $y^2=x^3+a_{4,E_0}(u)x+a_{6,E_0}(u)$
have the same
$j$-invariant, and therefore, there exists $v$ in the algebraic closure
of the fraction field of $R[[q]]$ such that $a=v^4a_{4,E_0}(u)$ and
$b=v^6a_{6,E_0}(u)$, respectively. From this point on, the rest of the argument
is as in the proof of Lemma \ref{thcan}.

Injectivity of $\iota_{E_0}$ is proved along the same lines.
\qed

Let us set
$$
\bM(R[[q]])_{good}=\{(a,b) \in \bM_{ord}(R[[q]]): \, a,b \in
R[[q]]^{\times}\}\, .
$$
Then we may consider the set of {\it good reduction solutions}
$$
\cU_{good}=\{(a,b) \in \bM(R[[q]])_{good} :\, f(a,b)=0\}\, .
$$
Similarly, we may consider the set of solutions of $f$ corresponding
to $E_0$,
$$
\cU_{E_0}=\{(a,b) \in \bM(R[[q]])_{E_0}: \, f(a,b)=0\}\, .
$$
Notice that, by Lemma \ref{L1}, we have
\begin{equation}
\label{defk}
\cU_{E_0} \subset \bM(R[[q]])_{CL,E_0}\, .
\end{equation}
Clearly,
\begin{equation}
\label{ciine}
\cU_{good}=\bigcup_{j(E_0) \neq 0,1728,ss} \cU_{E_0}\, ,
\end{equation}
where $ss$ are the supersingular values.
In order to study the set $\cU_{good}$, it suffices to study the sets
$\cU_{E_0}$ for all $E_0$ with $j(E_0) \neq 0,1728,ss$.
This is the object of the next Theorem.

 Let us define, for $v \in R[[q]]^{\times}$ and $\kappa \in \bZ_+$,
$\kappa \not\in p\bZ$, $z \in \bZ_p$, $\alpha \in R$,
the following pairs:
$$
\begin{array}{rcll}
u_{\Gamma,E_0,v} & := & (v^4 a_{4,E_0}(0), v^6 a_{6,E_0}(0)) & \in
\bM_{ord}(R[[q]])\, , \vspace{1mm} \\
u^z_{\Gamma,E_0,v,\kappa,\alpha} & := & (v^4
a_{4,E_0}(u^z_{m,\kappa,\alpha}-1), v^6
a_{6,E_0}(u^z_{m,\kappa,\alpha}-1)) & \in
\bM_{ord}(R[[q]])\, .
\end{array}
$$

Notice that, by Lemma \ref{L2}, the maps
$$
\begin{array}{rcl}
 v & \mapsto & u_{\Gamma,E_0,v}\, ,\\
(\alpha, v) & \mapsto & u^z_{\Gamma,E_0,v,\kappa,\alpha}\, ,
\end{array}
$$
are injective. By Remark \ref{iepure}, $u_{\Gamma,E_0,v}$ does not depend
on our choice of a basis of the Tate module $T_p(E_0)$. Also,
multiplying a basis of the Tate module by an element $\lambda \in
\bZ_p^{\times}$ has the effect of replacing
$u_{\Gamma,E_0,v,\kappa,\alpha}$ by
$u_{\Gamma,E_0,v,\kappa,\lambda^2 \alpha}$.

\begin{theorem}
\label{tzutzu}
Consider the arithmetic heat equation given by
$$
f=f^1_q+\lambda f^1_{p,-2}+\epsilon f^{1\phi}_{p,-2}\, ,
$$
where
$\lambda \in R^{\times}$, $\epsilon \in pR$, and
$z:=\frac{\epsilon}{\lambda} \in p\bZ_p$.
Let $E_0$ be an ordinary elliptic curve over $k$ with
$j \neq 0,1728$, $\kappa:=\frac{c_{E_0}}{\gamma_{E_0}} \cdot
\lambda \in R^{\times}$, and $\cU_{E_0}$ be the set of solutions
corresponding to $E_0$. Then the following hold:
\begin{enumerate}
\item Assume $\kappa  \not\in \bZ_+$. Then
$$
\cU_{E_0}=\{u_{\Gamma ,E_0,v} |\, v \in R[[q]]^{\times}\}\, .
$$
\item Assume $\kappa \in \bZ_+$. Then
$$
\cU_{E_0}=\{u^z_{\Gamma,E_0,v,\kappa,\alpha} |\, v \in
R[[q]]^{\times}, \, \alpha \in R\}\, .
$$
\end{enumerate}
\end{theorem}

\begin{remark}
According to the two cases above, we have corresponding canonical
identifications $\cU_{E_0}/R[[q]]^{\times} \simeq 0$ and
$\cU_{E_0}/R[[q]]^{\times} \simeq R$, where the latter is obtained via
the $\alpha$ parameter.
This and those identifications in (\ref{ggica})
show, for instance, that there are canonical identifications between bad
and good reduction solutions
$$
\cU_{\zeta}/R[[q]]^{\times} \simeq \cU_{E_0}/R[[q]]^{\times}
$$
of the arithmetic convection equations defined by
$$
f^1_q+\xi\frac{\gamma}{c} f^1_{p,-2}\quad \text{and}
\quad f^1_q+\xi\frac{\gamma_{E_0}}{c_{E_0}} f^1_{p,-2}\, ,
$$
respectively, where $\xi$ is any element of $R^{\times}$. \qed
\end{remark}

{\it Proof of Theorem {\rm \ref{tzutzu}}}.
By Lemmas \ref{L0}, \ref{L1}, \ref{L2}, and (\ref{defk}),
$\cU_{E_0}$ consists of all pairs $(v^4a_{4,E_0}(u), v^6
a_{6,E_0}(u))$ with $v \in R[[q]]^{\times}$, $u \in qR[[q]]$, such
that
\begin{equation}
\label{cococoo}
\Psi_q(1+u)+ \kappa \Psi_p(1+u)+\kappa z \phi \Psi_p(1+u)=0\, .
\end{equation}
The conclusion follows from this exactly as in the case of bad reduction.
\qed

\section{Instability}
\setcounter{theorem}{0}
We recall the series in (\ref{ursu}), and
start with the following obvious remark:
for any fixed $0 \neq \kappa \in \bZ_+$,  the function
$$ \begin{array}{rcl}
\bZ_p \times R  & \ra & R[[q]]\\
(z,\alpha)  & \mapsto & u^z_{a,\kappa,\alpha}
\end{array}$$
is continuous for the $p$-adic topologies
of $\bZ_p \times R$ and $R[[q]]$.
Morally speaking, this is an {\it stability} result for $\bG_a$.

In deep contrast with this, we have the following {\it instability} result
for $\bG_m$ and $\bM_{\Gamma}$:

\begin{theorem}
\label{instability} Let $\kappa \in \bZ_+$, $\kappa \not\in p\bZ$,
$\alpha_0, \alpha \in R$,
$\eta_0,\eta \in R^{\times}$,
 $z_0, z \in \bZ_p$.
Assume that one of the following holds:
\begin{enumerate}
\item $\alpha =\alpha_0$ and $z \neq z_0$,
\item $\alpha \neq \alpha_0$ and $z=z_0$.
\end{enumerate}
\noindent Then we have that
\begin{equation}
\label{non1}
\eta u^z_{m,\kappa,\alpha} \not\equiv \eta_0 u^{z_0}_{m,\kappa,\alpha_0}
\; \text{{\rm mod} $p$ in $R[[q]]$}\, .
\end{equation}
Moreover, if $v_0,v \in R[[q]]^{\times}$, then
\begin{equation}
\label{non2}
u^z_{\Gamma, \eta,v, \kappa,\alpha} \not\equiv
u^{z_0}_{\Gamma, \eta_0, v_0,\kappa,\alpha_0}\; \text{{\rm mod} $p$ in
$R[[q]]$}\, ,
\end{equation}
\begin{equation}
\label{muiee}
u^z_{\Gamma, E_0, v, \kappa,\alpha} \not\equiv
u^{z_0}_{\Gamma, E_0, v_0,\kappa,\alpha_0}\; \text{{\rm mod} $p$ in
$R[[q]]$}\, .
\end{equation}
\end{theorem}

Under the hypothesis (1), (\ref{non2}) morally implies
that for $z_0=0$, for instance, if $\alpha$
is a given ``boundary condition,'' then no matter how close the
``parameter'' $z\neq 0$ is to $0$, a solution
$u^z_{\Gamma, \zeta \beta,v, \kappa,\alpha}$ of the arithmetic heat equation
defined by $f^1_q+\lambda f^1_{p,-2}+\epsilon f^{1 \phi}_{p,-2}$
($z=\epsilon/\lambda$)
with ``boundary condition'' $\alpha$ is never close to a solution
$u^0_{\Gamma, \zeta_0 \beta_0, v_0,\kappa,\alpha}$ of the arithmetic
convection equation $f^1_q+\lambda f^1_{p,-2}$
that has the same ``boundary condition'' $\alpha_0=\alpha$.
A similar interpretation holds for (\ref{muiee}).
Furthermore, under the hypothesis (2), (\ref{non2}) implies that given
an arithmetic heat (or convection)
equation, any two solutions corresponding to different ``boundary
conditions'' $\alpha \neq \alpha_0$ are always far apart.

{\it Proof.} We assume hypothesis (1). The case when using (2) instead
can be treated in a similar way.

We take $\frac{\partial}{\partial x}$ in the recurrence relation
(\ref{recurrence}) for the $b_n(x,y)$s, and set $x=y=0$. We then
get
$$\frac{\partial b_n}{\partial x} (0,0)=\frac{\partial b_{n-1}}
{\partial x} (0,0) \; \text{for $n \geq 2$}\, .$$
Since $\frac{\partial b_1}{\partial x} (0,0)=-1$, we obtain
that
$$\frac{\partial b_n}{\partial x} (0,0)=-1 \;
\text{for $n \geq 1$}\, .$$
In particular, viewing the $b_n$s as elements in $\bZ[[x,y]]$,
for all $n\geq 1$, we may write
$$b_n(x,y)=b_n(0,y)-x+A_n(x,y)xy+B_n(x,y)x^2\, ,$$
where $A_n,B_n \in \bZ[[x,y]]$.
We deduce that
$$\begin{array}{rclll}
b_n(pz,p)-b_n(pz_0,p) & = & p(z-z_0) & + &
p^2z[A_n(pz,p)-A_n(pz_0,p)]\vspace{1mm} \\
\  & \  & \  & + & p^2(z-z_0)A_n(pz_0,p) \vspace{1mm} \\
\  & \  & \  & + & p^2z^2[B_n(pz,p)-B_n(pz_0,p)] \vspace{1mm} \\
\  & \  & \  & + & p^2(z^2-z^2_0)B_n(pz_0,p)\, .
\end{array}$$
Let us denote by $v_p:\bQ_p^{\times} \ra \bZ$ the $p$-adic valuation.
Since
$$\begin{array}{rcl}
v_p(A_n(pz,p)-A_n(pz_0,p)) & \geq  & v_p(z-z_0)+1\, , \\
v_p(B_n(pz,p)-B_n(pz_0,p)) & \geq  & v_p(z-z_0)+1\, ,
\end{array}
$$
we deduce that
\begin{equation}
\label{esti}
v_p(b_n(pz,p)-b_n(pz_0,p))=v_p(z)+1 \; \text{for all $n\geq 1$}\, .
\end{equation}

Assume now that (\ref{non1}) is false, that is to say,
$$\eta u^z_{m,\kappa,\alpha}-\eta_0 u^{z_0}_{m,\kappa,\alpha} \in pR[[q]]\, .$$
We derive a contradiction.
Setting $q=0$ in this expression, we get that
$$\eta \equiv \eta_0\; \text{mod $p$}\, ,$$
and it follows that
$$\frac{u^z_{m,\kappa,\alpha}}{u^{z_0}_{m,\kappa,\alpha}} \in 1+pR[[q]]\, .$$
On the other hand,
$$\frac{u^z_{m,\kappa,\alpha}}{u^{z_0}_{m,\kappa,\alpha}}=\exp\left(
\sum_{n \geq 0}(b_n(pz,p)-b_n(pz_0,p))\alpha^{\phi^n} \frac{q^{\kappa p^n}}
{\kappa p^n}\right)\, .$$
Using the fact that the function $X \mapsto \log X$ maps
$1+pR[[q]]$ into $pR[[q]]$,
we deduce that
$$\sum_{n \geq 0}(b_n(pz,p)-b_n(pz_0,p))\alpha^{\phi^n} \frac{q^{\kappa p^n}}
{\kappa p^n}\in pR[[q]]\, .$$
In particular, for any $n \geq 1$, we have that
$$v_p\left( (b_n(pz,p)-b_n(pz_0,p))\frac{\alpha^{\phi^n}}{\kappa p^n} \right)
 \geq 1\, ,$$
hence, by (\ref{esti}), that
$$v_p(z)+1+v_p(\alpha)-n \geq 1$$
for all $n \geq 1$. This is a contradiction.

We now prove (\ref{non2}) arguing again by contradiction. Let us assume
that (\ref{non2}) is false, that is to say, that we have that
$$
u^z_{\Gamma, \eta,v, \kappa,\alpha} \equiv
u^{z_0}_{\Gamma, \eta_0, v_0,\kappa,\alpha}\; \text{mod $p$ in $R[[q]]$}\, .
$$
Let $\frac{1}{j_{\infty}}(q) \in \bZ[[q]]$
  be the reciprocal of the Fourier expansion of the $j$-function.
Then
 $$\frac{1}{j_{\infty}}(\eta q u^z_{m,\kappa,\alpha}) \equiv
 \frac{1}{j_{\infty}}(\eta_0 q u^{z_0}_{m,\kappa,\alpha})\;
\text{mod $p$ in $R[[q]]$}\, .
$$
 Now $\frac{1}{j_{\infty}}(q)$ has a compositional inverse in $qR[[q]]$.
We conclude that
 $$
 \eta u^z_{m\kappa,\alpha} \equiv \eta_0 u^{z_0}_{m,\kappa,\alpha}\;
\text{mod $p$}\, ,
 $$
which contradicts (\ref{non1}). Thus, the assumption is false, and
(\ref{non2}) holds.

Another argument similar to the one above can be used to prove (\ref{muiee}).
\qed

\begin{remark} Theorem \ref{instability} implies in particular that
for any fixed
$\kappa \in \bZ_+$, $\kappa \not\in p\bZ$,  $z \in \bZ_p$, $v \in R[[q]]$,
$\alpha_0 \in R$ the functions
\begin{equation}
\label{2maps}
\begin{array}{rcl}
R  & \ra & R[[q]] \times R[[q]]\\
\alpha  & \mapsto & u^z_{\Gamma,\eta,v,\kappa,\alpha}\\
\end{array}\, ,
\quad
\begin{array}{rcl}
R  & \ra & R[[q]] \times R[[q]]\\
\alpha & \mapsto & u^z_{\Gamma,E_0,v,\kappa,\alpha}
\end{array}\, ,
\end{equation}
are discontinuous at $\alpha_0$ for the $p$-adic topologies
of $R$ and $R[[q]]\times R[[q]]$. However,
if $$c_n:R[[q]] \times R[[q]] \ra R\times R$$
($n \geq 1$) are the maps
$$c_n(\sum a_iq^i,\sum b_iq^i)=
(a_n,b_n)\, ,$$ then
 the maps
\begin{equation}
\label{puzza} \begin{array}{rcl}
R  & \ra & R \times R\\
\alpha  & \mapsto & c_n(u^z_{\Gamma,\eta,v,\kappa,\alpha})\\
\end{array}\, ,
\quad
\begin{array}{rcl}
R  & \ra & R \times R\\
\alpha & \mapsto & c_n(u^z_{\Gamma,E_0,v,\kappa,\alpha})
\end{array}\, ,
\end{equation}
are continuous for the $p$-adic topologies.
Indeed,
the maps (\ref{2maps}) are pseudo $\delta_p$-polynomial in the sense of the
Definition 6.10 in \cite{pde}. In particular, we see
that the family (indexed by $n \geq 1$) of maps (\ref{puzza}) is {\it not}
equicontinuous at any $\alpha_0$.
\end{remark}

\section{Analogies between various equations}

The following table summarizes some of the analogies between
various arithmetic convection equations introduced in \cite{pde} and in
the present paper. We freely employ the notation of \cite{pde},
and that introduced here.
\bigskip

\begin{center}
\begin{tabular}{||l|l|l|l||} \hline \hline
\text{variety} & \text{equation} & \text{stationary} &
\text{solutions with fixed}\\
\  & \  & \text{solutions}  & \text{behavior at infinity}\\
\hline
 $\bG_a $ & $\dd+\lambda
\phi_p-\lambda$ & $\cU_0=\bZ_p$ & $\cU_1\simeq R$ or $0$\\
\hline $\bG_m$ & $\psi_q+\lambda \psi_p$ & $\cU_0=\mu(R)$ &
$\cU_1\simeq R$ or $0$ \\
\hline $E={\rm Tate}(\beta q)$ & $\psi^1_{pq}$ & $\cU_0^1$ & 
$\cU_1^1\simeq R$ or
$0$ \\
\hline $E/R$, a CL curve& $\psi_q^1+\lambda \psi_p^1$ &
$\cU_0=\cap p^n E(R)$ & $\cU_1\simeq R$ or $0$ \\
\hline $\bM_{\Gamma}$  & $f^1_q+\lambda
f^1_{p,-2}$ &
$\cU_0=\emptyset$ & $\cU_{\zeta}/R[[q]]^{\times} \simeq R$ or $0$\\
\  & \  & (bad reduction) & (bad reduction)\\
 \  & \  & $\cU_0=\bM(R)_{CL}$ &
$\cU_{E_0}/R[[q]]^{\times} \simeq R$ or $0$ \\
\  & \  & (good reduction) & (good reduction)\\
\hline
\end{tabular}
\end{center}

\bigskip

In the table above, the set of stationary solutions $\cU_0$ for
$\bM_{\Gamma}$, in the bad reduction case, is defined as the
solutions in $\cU_{bad} \cap \bM_{ord}(R)$; the latter set is, of
course empty. Also, the set $\cU_0^1$ of {\it stationary solutions}
for $E={\rm Tate}(\beta q)$, defined in \cite{pde} consists of all
solutions in the Manin kernel; cf. \cite{pde} for details. The
alternative ``$R$ or $0$'' in the last column depends on the value of
$\lambda$ (for all cases except $E={\rm Tate}(\beta q)$) respectively
$\beta$ (in the case of $E={\rm Tate}(\beta q)$).

 A similar table can
be provided for arithmetic heat equations of order $2$
(in which $E/R$ is, this time, non-CL). On the
other hand, as was shown in this paper, the analogy breaks down for
wave equations: these wave equations exist (and have interesting
 solutions) for $\bG_a,\bG_m$, and $E$, but do not exist for $\bM_{\Gamma}$.

\bibliographystyle{amsalpha}

\end{document}